\pdfoutput=1
\documentclass[a4paper, 10pt]{article}

\usepackage{preamble}
\usepackage{mathtools}
\DeclareMathOperator{\Null}{Null}

\title{A Theorem of the Heart for the K-theory of Endomorphisms}
\date{\today}
\author{Victor Saunier}

\begin{document}
\maketitle

\begin{abstract}
    We show that Quillen's resolution theorem for K-theory also applies to exact $\infty$-categories. We introduce heart structures on a stable $\infty$-category, generalizing weight structures, and using resolution ideas, we show that the category of stable $\infty$-categories equipped with a heart structure fully-faithfully embeds into the category of exact $\infty$-categories. Consequently, we show a generalized theorem of the heart for K-theory, which is equivalent to its invariance under passage to the stable envelope of exact $\infty$-category in the image of the heart functor. Finally, leveraging the above, we show that K-theory of endomorphisms satisfies the theorem of the heart for weight structures, even allowing coefficients in a suitable bimodule.
\end{abstract}

\tableofcontents

\vspace{2.5em} 

The resolution theorem is one of the three major results proven by Quillen \cite[Theorem 3]{Quillen} when he introduced the higher K-theory of exact categories. It states the following: let $\Ccal$ be an exact 1-category and $\Acal$ an exact subcategory closed under extensions in $\Ccal$ satisfying the following two conditions: 
\begin{itemize}
    \item[(i)] For every exact sequence $\begin{tikzcd}[cramped]X\arrow[r, hook] & Y\arrow[r, two heads] & Z\end{tikzcd}$ with $Y\in\Acal$, we have $X\in\Acal$
    \item[(ii)] For every $Z\in\Ccal$, there is an exact sequence $\begin{tikzcd}[cramped]X\arrow[r, hook] & Y\arrow[r, two heads] & Z\end{tikzcd}$ with $Y\in\Acal$.
\end{itemize}
then, the inclusion $\Acal\to\Ccal$ induces an equivalence $\Kth(\Acal)\xrightarrow{\simeq}\Kth(\Ccal)$. The above axioms in particular imply that every object of $\Ccal$ admits a resolution of length 1 by objects of $\Acal$, hence the name. In Corollary 1 of \textit{loc. cit.}, Quillen extends the result to the situation where objects of $\Ccal$ have finite-length resolutions by objects of $\Acal$, under similar suitable hypotheses on the inclusion. 

The standard example given by Quillen of a \textit{resolving} situation is the inclusion of the category of compact projective $R$-modules inside the category of compact $R$-modules for some ring $R$. \\

Compared to its companion facts, the localization and devissage theorems which Quillen also proved in \textit{loc. cit.}, the resolution theorem has since enjoyed less time under the spotlight. 

In \cite[1.4]{Waldhausen}, Waldhausen provided a version of localization adapted to Waldhausen 1-categories; he called it the additivity theorem and under this name, the localization property has grown to become the cornerstone of so-called additive and localizing invariants, furnishing a universal property for K-theory of stable categories following Blumberg-Gepner-Tabuada in \cite{Blumberg}, allowing the extension of the definition to large dualizable categories by Efimov, as well as those ideas applying to L-theory and Grothendieck-Witt theory in spectacular fashion (see the series of papers \cite{HermKI, HermKII, HermKIII, HermKIV}). 

As for devissage, Waldhausen's result that looks the most like it is the cell filtration theorem, also known as Waldhausen's sphere theorem \cite[Theorem 1.7.1]{Waldhausen}. This is not a generalization, as noted by Thomason-Trobaugh in \cite[1.11.1]{ThomasonTrobaugh}, but a result with a similar flavour. A celebrated theorem of Gillet-Waldhausen (see Theorem 1.11.7 \textit{loc. cit.}) is often claimed as some ersatz of devissage, comparing the K-theory of an exact 1-category $\Ecal$ nicely embedded into an abelian one with the 1-category of bounded chain complexes in $\Ecal$, obtained by restricting those in the ambient abelian one. More recently, Raptis gave a generalization of Quillen's devissage and other devissage-like statements in the 1-categorical world \cite[Theorem 5.7, Theorem 6.11]{Raptis}, the latter of which is general enough to recover higher categorical ones. \\

In the higher categorical world, Barwick showed in \cite{BarwickHeart} that if $\Ccal$ is a stable category with a bounded t-structure, then the inclusion of its heart $\Ccal^\heart\to\Ccal$ is sent to an equivalence by K-theory. Barwick named this result the theorem of the heart, in reference to an earlier result of Neeman for triangulated categories \cite{NeemanHeartIIIA, NeemanHeartIIIB}.

The central example of a theorem-of-the-heart situation is obtained by realizing any abelian 1-category as the heart of a t-structure on its bounded derived category. Note that although the statement looks similar, the major difference to the Gillet-Waldhausen setting is the apparition of $\infty$-categories. Like Gillet-Waldhausen however, Barwick's proof relies on his generalized version of a localization theorem. Let us also mention that in the context of stable categories and t-structure, Burklund and Levy have shown a devissage-like result in \cite{BurklundLevy}, by combining the better of both Barwick's theorem of the heart and Quillen's original devissage result. \\

This leaves Quillen's resolution theorem lacking a major counterpart in the setting of stable categories. Coming to the rescue is another theorem of the heart, whose first correct proof can be found in work of Hebestreit-Steimle \cite[Corollary 8.1.3]{HebestreitSteimleAppHarpaz} where is proven the hermitian version of this fact. 

This theorem has the same statement as Barwick's, but concerns weight structures instead of t-structures: if $\Ccal$ is a stable category equipped with a bounded weight structure, then the inclusion of its heart $\Ccal^\heart\to\Ccal$ is sent to an equivalence by K-theory, where the heart is an additive category endowed with the split-exact structure. Unlike Barwick's, this theorem of the heart applies to the inclusion of the category of compact projective modules into the bounded derived category, which serves as the higher categorical version of the category of compact modules over some ring $R$. This recovers Quillen's initial example thanks to Barwick's theorem of the heart. More generally, this  applies to any additive 1-category into its category of bounded chain complex by \cite[Example 3.1.3]{HebestreitSteimleAppHarpaz}. \\

This article holds the following belief to be true:

\begin{slogan}
    The theorem of the heart of \cite{HebestreitSteimleAppHarpaz} is the correct higher analogue of Quillen's resolution theorem.
\end{slogan}

In this slogan, we mean both the theorem of the heart of weight structures of Hebestreit-Steimle and the generalization we will introduce. We will give an explicit version of this slogan in Theorem \ref{HeartIsResolutionIntro} (see \ref{HeartIsResolution} in the text). Let us explain our strategy. \\

By a result of Sosnilo \cite[Corollary 3.4]{SosniloHeart}, the heart is a fully-faithful functor from stable categories equipped with a bounded weight-structure into the category of additive categories. Its essential image consists of weakly-idempotent additive categories and for such categories, the inverse is given by mapping an additive $\Acal$ to its stable envelope $\Stab(\Acal)$ endowed with a suitable weight structure. 

Additive categories can be considered as exact categories via their closed-under-extension embedding into their stable envelope, but as the above suggest, if $\Acal$ is the heart of a bounded weight structure on a stable category, every exact sequence of $\Acal$ is split so this structure is the same as the split-exact one. This is unfortunate as Quillen's resolution theorem lives quite naturally in the world of exact categories, where all sequences need not split, and if one is convinced by our slogan, then one will agree that a piece of the puzzle is missing. \\

We claim this missing piece is given the following definition: 

\begin{defi}
    Let $\Ccal$ be a stable category. A \textit{heart structure} on $\Ccal$ is the datum of a pair of full subcategories $(\Ccal_{\geq 0}, \Ccal_{\leq 0})$ subject to the following conditions:
    \begin{itemize}
        \item{(i)} $\Ccal_{\geq 0}, \Ccal_{\leq 0}$ are closed under extensions, $\Ccal_{\geq 0}$ under finite colimits and $\Ccal_{\leq 0}$ under finite limits.
        \item{(iii)\footnote{This indexing is intentional.}} For any $C\in\Ccal$, there is an exact sequence 
        $$
            \begin{tikzcd}[cramped]
                X\arrow[r] & C\arrow[r] & \Sigma Y
            \end{tikzcd}
        $$
        with $X\in\Ccal_{\leq 0}$ and $Y\in\Ccal_{\geq 0}$.
    \end{itemize}
    
    A \textit{heart-exact} functor between heart structures is an exact functor $f:\Ccal\to\Dcal$ on the underlying categories which preserves the two full subcategories of the structure. We denote $\CatExHeart$ the category of heart structures and heart-exact functors between them. 
\end{defi}

The aforementioned weight structures are in particular heart structures, as they have the stronger property of satisfying an axiom (ii), guaranteeing that the mapping spectra $\map_\Ccal(X, Y)$ are connective for $X\in\Ccal_{\leq 0}$ and $Y\in\Ccal_{\geq 0}$, as well as connective and coconnective objects being closed under retracts. 

It still makes sense to talk about boundedness for heart structures, just as it still makes sense to talk about their heart (thankfully!). However, it is no longer natural to consider the heart of a heart category only as an additive category. Instead, we remark that $\Ccal^\heart$ is closed under extensions in $\Ccal$, hence inherits an exact structure given by specifying a sequence to be exact in $\Acal$ if and only if it is exact in $\Ccal$. This is a sensible idea: by removing axiom (ii) of the definition of a weight-structure, we are no longer guaranteed that every exact sequence of the heart splits. In particular, a stable category with a heart structure need not be recovered by the additive heart. However, considering the heart as an exact category will suffice, which is the point of our first result:

\begin{thm}
    The functor $(-)^\heart:\CatExHeartBounded\to\ExactCat$ taking a bounded heart category to its heart is fully-faithful. 
\end{thm}

For abstract reasons, the inclusion $\CatEx\to\ExactCat$ admits a left adjoint, the \textit{stable envelope} $\Stab(\Ecal)$. This stable envelope coincides with the stable envelope of the underlying additive categories when the structure is split-exact. By \cite[Theorem 1]{Klemenc}, the canonical functor $\Ecal\to\Stab(\Ecal)$ is fully-faithful, reflects exact sequences and its image is closed under extensions. 

The above theorem follows from showing that if $\Ccal$ has a bounded heart structure, then, the inclusion of its heart $\Ccal^\heart\to\Ccal$ induces an equivalence $\Stab(\Ccal^\heart)\simeq\Ccal$, i.e. $\Ccal^\heart$ generates $\Ccal$ as an exact category. To prove this, we show an explicit criterion for an exact functor to induce an equivalence of stable envelopes; interestingly, the key input is provided by the hypotheses of Quillen's resolution theorem:

\begin{lmm}
    Suppose $i:\Acal\to\Ccal$ is a \textit{resolving} functor, i.e. satisfies the conditions of Quillen's Theorem 3 in \cite{Quillen}. Then, for every stable $\Dcal$, the following functor
    $$
        \begin{tikzcd}
            i^*:\FunEx(\Ccal, \Dcal)\arrow[r, "\simeq"] & \FunEx(\Acal, \Dcal)
        \end{tikzcd}
    $$
    is an equivalence. In particular, $i$ induces an equivalence $\Stab(\Acal)\xrightarrow{\simeq}\Stab(\Ccal)$.
\end{lmm}

One can show that if $\Ccal$ has a heart structure, then $\Ccal^\heart\to\Ccal$ is a composition of resolving and op-resolving\footnote{Functors $i$ such that $\catop{i}$ is resolving} functors, possibly infinite. Since $\Ccal$ is stable, the above lemma implies that the map $\Stab(\Ccal^\heart)\simeq\Ccal$ is an equivalence. \\

Unlike Sosnilo's result, we do not know of a completion description of the essential image of our heart functor. Hence, let us introduce a bit of terminology:

\begin{defi}
    An exact category $\Ecal$ is said to be \textit{stably-exact} if it is the heart of a heart category.
\end{defi}

Additive categories are always stably-exact, by Sosnilo's result in \cite[Corollary 3.4]{SosniloHeart}, and stable categories are also stably-exact when considered as exact categories, since they are equivalent to their stable envelope and thus one can take the maximal structure where every object is both connective and coconnective. We are also able to show that if $\Acal$ is the heart of a heart structure on $\Ccal:=\Stab(\Acal)$, then all the intermediary categories $\Ccal_{[n, m]}:=\Ccal_{\leq m}\cap\Ccal_{\geq n}$ for integers $n\leq 0\leq m$ are also stably-exact, providing many non-trivial examples. \\

As a consequence of the proof of the above theorem, we get the following conceptual, modern interpretation of the conditions of Quillen's resolution theorem, at least between stably-exact categories:

\begin{thm}[Resolution-Heart equivalence] \label{HeartIsResolutionIntro}
    Let $F:\ExactCat\to\Ecal$ be a functor preserving filtered colimits. The following are equivalent:
    \begin{itemize}
        \item[(i)] For every stably-exact $\Ecal$, $F$ sends the map $\Ecal\to\Stab(\Ecal)$ to an equivalence.
        \item[(ii)] For every bounded heart category $\Ccal$, $F$ sends the map $\Ccal^\heart\to\Ccal$ to an equivalence.
        \item[(iii)] $F$ sends every resolving and every op-resolving functor $\Acal\to\Dcal$ with $\Acal, \Dcal$ stably-exact to an equivalence.
    \end{itemize}
\end{thm}

Though this does not speak about resolving functors between non-stably exact categories, we believe this justifies the previously stated Slogan. \\

We are also able to show that Quillen's resolution theorem also holds for K-theory of exact $\infty$-categories, by an almost identical adaptation of Quillen's original arguments in \cite[Theorem 3, Corollary 1]{Quillen}. In particular, using previously mentioned ideas, we get a generalized theorem of the heart for K-theory, which applies to heart structures:

\begin{thm}[Theorem of the heart]
    Let $\Ccal$ be a bounded heart category and denote $\Ccal^\heart$ its heart. Then, the map
    $$
        \begin{tikzcd}
            \Kth(\Ccal^\heart)\arrow[r, "\simeq"] &\Kth(\Ccal)
        \end{tikzcd}
    $$
    is an equivalence, where on the left hand side, we have taken K-theory of the exact category $\Ccal^\heart$. Equivalently, for every stably-exact $\Ecal$, the following map is an equivalence:
    $$
        \begin{tikzcd}
            \Kth(\Ecal)\arrow[r, "\simeq"] &\Kth(\Stab(\Ecal))
        \end{tikzcd}
    $$
\end{thm}

The application we give of our generalized theorem of the heart was at the genesis of this whole article. If $\Ccal$ is a stable category and $M:\Ccal\to\Ind\Ccal$ a $\Ccal$-bimodule (i.e. an exact functor), we let $\Lace(\Ccal, M)$ denote the lax-equalizer from the Yoneda embedding $j:\Ccal\to\Ind\Ccal$ to $M$. This is a stable category.

If $\Ccal\simeq\Perf(R)$, the category of compact $R$-modules over some ring spectra $R$, and $M$ a $R$-bimodule seen as a functor to the $\Ind$-construction via tensoring $M\otimes_R-$, then $\Lace(\Ccal, M\otimes_R-)$ is exactly $\End(R, M)$, the category of $M$-parameterized endomorphisms. When $R$ and $M$ are further supposed connective, $\End(R, \Sigma M)$ coincides with $\Perf(R\oplus M)$, the category of compact modules over the square-zero extension (see Example \ref{SqZeroVSLace}). Remark also that if $M=j$ is the Yoneda embedding, then $\Lace(\Ccal, M)$ coincides with $\End(\Ccal)$. We will call the composite $\Kth(\Lace(\Ccal, M))$ \textit{laced K-theory}. \\

If $\Ccal$ has a heart structure, we say that a $\Ccal$-bimodule $M$ is \textit{weighted} if, for every $X\in\Ccal_{\leq 0}$ and $Y\in\Ccal_{\geq 0}$, the mapping spectra $\map(X, M(Y))$ is connective. When $\Ccal=\Perf(R)$ and the heart structure is the usual weight structure with heart $\Proj(R)$, weighted bimodules are exactly connective $R$-bimodules. 

When $M$ is weighted, $\Lace(\Ccal, M)$ inherits a heart structure whose heart is the full subcategory $\Lace(\Ccal^\heart, M)$ of pairs $(X, f)$ where $X\in\Ccal^\heart$. Note that even if $\Ccal$ has a heart structure, the best one can hope for in general is a heart structure on $\Lace(\Ccal, M)$ which verifies axiom (ii) of weight structures only up to a shift by $-1$, hence no longer a weight structure. Still, our generalized theorem of the heart yields, for every heart category $\Ccal$ and every weighted bimodule $M$, an equivalence:
$$
    \begin{tikzcd}
        \Kth(\Lace(\Ccal^\heart, M))\arrow[r, "\simeq"] & \Kth(\Lace(\Ccal, M))
    \end{tikzcd}
$$
Note however that it need not be that if $\Ccal^\heart$ has a heart structure and $M$ is any bimodule, then $\Lace(\Ccal^\heart, M)$ is the heart of a heart structure on $\Lace(\Ccal, M)$. In fact, we expect the above equivalence to fail in general if $\Ccal$ has only a $(-1)$-shifted weight structure as above and $M$ is not weighted. \\

Remark that the Yoneda embedding is a weighted bimodule if and only if the heart structure on $\Ccal$ verifies the condition (ii) in the definition of a weight structure. Hence we have:

\begin{thm}[Theorem of the heart for $\Kth^{\End}$] \label{KEndHeart}
    Let $\Ccal$ be a stable category equipped with a bounded weight structure. Then, there is an equivalence
    $$
        \begin{tikzcd}
            \Kth(\End(\Ccal^\heart))\arrow[r, "\simeq"] & \Kth(\End(\Ccal))
        \end{tikzcd}
    $$
\end{thm}

Equivalently, the above states that if $\Acal$ is additive, then $\Kth(\End(\Acal))\to\Kth(\End(\Stab(\Acal)))$ is an equivalence, though one has to be careful and $\End(\Acal)$ is not considered simply as an additive category but as an exact category with the inherited structure from $\End(\Stab(\Acal))$. This generalizes Theorem 1.7 of \cite{BlumbergGepnerTabuadaKEnd}, which proved the above for $\Acal$ being the 1-category of projective modules over a discrete ring $A$. It also recovers Theorem 4.9 of \textit{loc. cit.} which showed the above held on $\pi_0$ for the category of projective modules over a connective ring spectrum. \\

There is another consequence worth mentioning of the theorem of the heart for laced K-theory. Using this result, we can show the functor $M\mapsto\Kth(\Lace(\Ccal, M))$ commutes with sifted colimits when restricted to weighted bimodules $M$. A similar statement plays usually a critical role in the Dundas-Goodwillie-McCarthy, and we expect that to use this generalized statement in the upcoming article revisiting trace methods.

\vspace{1.2em} 

\textbf{Acknowledgements.} We want to thank Yonatan Harpaz for his constant support during the writing of this article, as well as helpful comments and insightful discussions. 

\vspace{1.2em} 

\textbf{Conventions.} In this article and as we have already done in this introduction, we write category for $\infty$-categories and specify 1-category when talking about categories whose mapping spaces are discrete. We write $\Spaces$ for the category of spaces, $\Sp$ for the category of spectra, $\CatInfty$ for the category of small categories, $\CatEx$ for its non-full subcategory of stable categories and exact functors between them and $\ExactCat$ for the category of exact categories and exact functors between them. If $\Ccal$ is a small category, we denote $\Ind\Ccal$ its $\Ind$-construction; in particular, if $\Ccal$ is stable, so is its $\Ind$-construction and we have $\Ind\Ccal\simeq\FunEx(\catop{\Ccal}, \Sp)$. When $\Ccal$ is stable, we denote by $\map$ the enrichment in $\Sp$ of the mapping spectra and reserve $\Map$ to mapping spaces.

\section{Quillen's resolution theorem for exact categories}

\hspace{1.2em} The goal of this section is to prove a higher categorical version of the original resolution theorem of Quillen, namely Theorem 3 of \cite{Quillen}. We follow a similar strategy, relying on the higher categorical version of Quillen's Theorem A provided by \cite[Theorem 4.1.3.1]{HTT} as well as the treatement of higher exact categories of \cite{BarwickHeart}. \\

\begin{defi}[Definition 3.1 of \cite{BarwickHeart}] \label{ExactCatDef}
    Let $\Ccal$ be an additive category. An \textit{exact structure} on $\Ccal$ is the datum of two subcategories $\Ccal_{inj}$ and $\Ccal_{proj}$ containing all equivalences and whose arrows are denoted respectively $\hookrightarrow$ and $\twoheadrightarrow$, which we choose to call respectively \textit{exact inclusions} and \textit{exact projections}, subject to the following conditions:
    \begin{itemize}
        \item For any $X\in\Ccal$, $0\hookrightarrow X$ is an exact inclusion and $X\twoheadrightarrow 0$ an exact projection.
        \item Exact inclusions are stable under pushout against any map and exact projections under pullback against any map\footnote{In particular, such pullbacks and pushouts are required to exist.}.
        \item Any square:
        $$
            \begin{tikzcd}
                X\arrow[r]\arrow[d] & Y\arrow[d] \\
                Z\arrow[r] & T
            \end{tikzcd}
        $$
        is a pullback of a span with a leg in $\Ccal_{inj}$ and the other in $\Ccal_{proj}$ if and only if it is a pushout of a cospan with the same condition.
    \end{itemize}
    
\end{defi}

If $\Ccal$ is an exact structure on an additive category, we will often omit to write the two subcategories and simply say that $\Ccal$ is an exact category. A product-preserving functor $\Ccal\to\Dcal$ between exact categories is itself called exact if it preserves either exact inclusions and their pushouts or exact projections and their pullbacks, in which case it preserves both by \cite[Proposition 4.8]{BarwickHeart}. \\

If $\Ccal$ is an exact category, an exact sequence in $\Ccal$ is a fiber-cofiber sequence of the form 
$$
    \begin{tikzcd}
        X\arrow[r, hook] & Y\arrow[r, two heads] & Z
    \end{tikzcd}
$$
Remark that it suffices that such a sequence is either a fiber or a cofiber sequence for the other one to also be satisfied. \\

Let $\Acal$ be an additive full subcategory of some exact category $\Ccal$ and suppose further that $\Acal$ is closed under extension in $\Ccal$. Then, there is an exact structure on $\Acal$ where exact inclusions are morphisms which are exact inclusions when viewed in $\Ccal$ and whose cofiber lies in $\Acal$, and projections are defined dually (this is marginally more general than 1.4.6 of \cite{BarwickQConstr}).

\begin{defi}
    Let $\Acal$ be an exact subcategory of $\Ccal$ closed under extensions. We say that the exact structure on $\Acal$ is \textit{inherited from $\Ccal$} if it is obtained by the above procedure. Specifically, a sequence is exact in $\Acal$ if and only if it is exact when viewed as a sequence in $\Ccal$.
\end{defi}

The inherited structure is the maximal exact structure such that $\Acal\to\Ccal$ is exact. The minimal such structure is the split-exact structure where a sequence is exact in $\Acal$ if and only if it is a split-exact sequence. \\

Associated to any exact category $\Ccal$ is another category denoted $Q(\Ccal)$, obtained by taking the underlying category of a complete Segal space $Q_\bullet(\Ccal)$ (see \cite[Section 3]{BarwickQConstr} or \cite[Section 2.7]{HermKII} for details). This category has the same objects and a morphism from $X$ to $Y$ is the datum of a span of the following shape:
$$
    \begin{tikzcd}[cramped]X & Z\arrow[r, hook]\arrow[l, two heads] & Y\end{tikzcd}
$$
The K-theory space is obtained from $Q(\Ccal)$ by looping once its geometric realization (i.e. the left adjoint of the inclusion $\Spaces\to\CatInfty$). In particular, an equivalence at the level of the Q-construction implies an equivalence at the level of K-theory. \\

In \cite[Theorem 3]{Quillen}, Quillen gives an explicit criterion for the embedding of a closed-under-extension $\Acal$ inside $\Ccal$ to induce an equivalence at the level of the Q-construction. We put it in a definition:

\begin{defi} \label{ResolvingDef}
    Let $i:\Acal\to\Ccal$ be a fully-faithful functor whose image is closed under extensions. Suppose $\Ccal$ is exact and endow $\Acal$ with the inherited structure. We say that $i$ is \textit{resolving} if it satisfies the following two additional properties:
    \begin{itemize}
        \item[(i)] For every exact sequence $\begin{tikzcd}[cramped]X\arrow[r, hook] & Y\arrow[r, two heads] & Z\end{tikzcd}$ with $Y\in\Acal$, we have $X\in\Acal$
        \item[(ii)] For every $Z\in\Ccal$, there is an exact sequence $\begin{tikzcd}[cramped]X\arrow[r, hook] & Y\arrow[r, two heads] & Z\end{tikzcd}$ with $Y\in\Acal$.
    \end{itemize}
    We say $i$ is op-resolving if $\catop{i}$ is resolving.
\end{defi}

Our naming convention is inspired by \cite{DundasTheorems} but beware that they call resolving an embedding satisfying conditions related to Quillen's Corollary 1, whereas we gave a name to the conditions of Quillen's Theorem 3. 

\begin{thm}[Quillen's resolution theorem] \label{QuillenTheoremThree}
    Let $i:\Acal\to\Ccal$ be a resolving functor. Then, $i$ induces an equivalence of spaces $|Q(\Acal)|\xrightarrow{\simeq}|Q(\Ccal)|$ so consequently an equivalence on the level of K-theory.
\end{thm}
\begin{proof}
    Denote $\Bcal$ the full subcategory of $Q(\Ccal)$ spanned by the image of $Q(\Acal)$. We have a factorization
    $$
        \begin{tikzcd}
            Q(\Acal)\arrow[r, "g"] & \Bcal\arrow[r, "f"] & Q(\Ccal)
        \end{tikzcd}
    $$
    and we will show that both maps induce equivalences on the geometric realizations by showing that $\catop{g}$ and $f$ are cofinal. For both proofs, we use the higher categorical version of Quillen's Theorem A (see \cite[Theorem 4.1.3.1]{HTT}), and are reduced to show some categories are weakly contractible. \\

    We first show $f$ is a weak equivalence. Let $X\in Q(\Ccal)$, it suffices to show that $\Mcal:=\Bcal\times_{Q(\Ccal)}Q(\Ccal)_{X/}$ is contractible; by construction, $\Mcal$ is a category whose objects are spans $X\twoheadleftarrow Z\hookrightarrow A$ with $A\in\Acal$. Consider the wide subcategory $Q^{proj}(\Ccal)$ of $Q(\Ccal)$ composed of spans $X\twoheadleftarrow Z\hookrightarrow Y$ where the map $Z\to Y$ is an equivalence; this category is equivalent to the wide subcategory $\catop{(\Ccal^{proj})}$ of $\catop{\Ccal}$ whose arrows are exact projections. 
    
    Denote $\Mcal^{proj}$ the subcategory of $\Mcal$ given by $\Bcal\times_{Q(\Ccal)}Q^{proj}(\Ccal)_{X/}$; this a full subcategory of $\Mcal$ since it is equivalently the pullback along the projection from $\Mcal$ of $Q^{proj}(\Ccal)_{X/}\to Q(\Ccal)_{X/}$ which is easily checked to be fully-faithful. Its objects are equivalently exact projections $A\twoheadrightarrow X$ with $A\in\Acal$ and a map from $A\twoheadrightarrow X$ to $A'\twoheadrightarrow X$ is the datum of an exact projection $A'\twoheadrightarrow A$ and a homotopy which makes the obvious triangle commute (note the order reversal).
    
    By \cite[Proposition 4.9]{HausgengHebestreitLinskensNuiten}, the collection of spans where the right hand map an equivalence (purely forward pointing spans in the parlance of \textit{loc. cit.}) and spans where the left hand map is the identity (purely backward pointing spans) forms an orthogonal factorization system on $Q(\Ccal)$ (see \cite[Definition 5.2.8.8]{HTT}). It follows from Lemma 5.2.8.19 of \textit{loc. cit.} that the inclusion $Q^{proj}(\Ccal)_{X/}\to Q(\Ccal)_{X/}$ admits a right adjoint given on objects as follows:
    $$
        \begin{tikzcd}
            (X & \arrow[l, two heads]Y\arrow[r, hook] & Z)\arrow[rr, mapsto] && (X & \arrow[l, two heads]Y\arrow[r, equal, "\id_Y"] & Y)
        \end{tikzcd}
    $$
    For any span $X\twoheadleftarrow Z\hookrightarrow Y$ with $Y\in\Acal$, we have $Z\in\Acal$ as well hence the above descends to a functor $\Mcal\to\Mcal^{proj}$. One readily checks that it provides a right adjoint to the inclusion $\Mcal^{proj}\to\Mcal$ hence the latter is a homotopy equivalence. Thus we are reduced to showing that $\Mcal^{proj}$ is weakly contractible.

    By (ii), the category $\Mcal^{proj}$ is nonempty since there exists $B\twoheadrightarrow X$ with $B\in\Acal$. Moreover, given two projections $Y\twoheadrightarrow X$ and $Y'\twoheadrightarrow X$, then closure under extensions of $\Acal$ implies that $Y\times_X Y'$ is an object of $\Acal$ since it fits in the following exact sequence of $\Ccal$:
    $$
        \begin{tikzcd}[cramped]
            Z\arrow[r] & Y\times_X Y'\arrow[r] & Y
        \end{tikzcd}
    $$
    where we have taken $Z$ to be the object fitting in an exact sequence $Z\hookrightarrow Y\twoheadrightarrow X$, which by (i) implies that $Z\in\Acal$. In particular, by fixing some object $p:A_0\twoheadrightarrow X$ in $\Mcal^{proj}$, and denoting $P$ the functor sending $A\twoheadrightarrow X$ to $A\times_{X}A_0\twoheadrightarrow X$, one obtains two natural transformations $\id\implies P$ and $\cst_{p}\implies P$ which together imply that $\Mcal^{proj}$ is weakly contractible. \\

    In order to show that $g$ is a weak equivalence, we let $X\in\Bcal$ and working dually, we show that $\Ncal:=Q(\Acal)\times_{\Bcal}\Bcal_{/X}$ is weakly contractible; this time, $\Ncal$ is a category whose objects are spans $Y\twoheadleftarrow Z\hookrightarrow X$ (mind the order) and all three $X, Y, Z\in\Acal$ though the maps of the span are only exact inclusions/projections in the exact structure on $\Ccal$. However, note that a map in $\Ncal$ from $(Y\twoheadleftarrow Z\hookrightarrow X)$ to $(Y'\twoheadleftarrow Z'\hookrightarrow X)$ are given by diagrams of the following shape:
    $$
        \begin{tikzcd}
            Y & T\arrow[l, two heads, "\in\Acal"]\arrow[d, hook, "\in\Acal"] & Z\arrow[l, two heads]\arrow[d, hook] \\
            & Y' & Z'\arrow[l, two heads]\arrow[d, hook] \\
            & & X
        \end{tikzcd}
    $$
    where every object is in $\Acal$, the maps indicated $\in\Acal$ are exact injections/projections in the exact structure of $\Acal$ and the square is cartesian in $\Ccal$; in particular, it follows from this last condition that $Z\hookrightarrow Z'$ is already an exact inclusion in the exact structure of $\Acal$. 
    
    Dually to what have done before, we can consider $Q^{inj}(\Ccal)$ the wide subcategory of $Q(\Ccal)$ whose arrows are spans $Y\twoheadleftarrow Z\hookrightarrow X$ with $Z\to Y$ an equivalence, and $\Bcal^{inj}$ the full subcategory of $Q^{inj}(\Ccal)$ spanned by $Q^{inj}(\Acal)$. Let $\Ncal^{inj}$ be the subcategory of $\Ncal$ given by $Q(\Acal)\times_{\Bcal}\Bcal^{inj}_{/X}$; again, $\Ncal^{inj}$ is a full subcategory of $\Ncal$. It has equivalently objects exact inclusions $Z\hookrightarrow X$ and morphisms given by exact inclusions $Z\xhookrightarrow{\in\Acal} Z'$ featuring in an exact sequence of $\Acal$ as well as a homotopy which makes the obvious diagram commute (note there is no order reversal here). 
    
    Finally, remark that the left adjoint to the inclusion $Q^{inj}(\Ccal)_{/X}\to Q(\Ccal)_{/X}$ which only keeps non-trivial the exact-inclusion, sends $\Bcal_{X/}$ to $\Bcal^{inj}_{X/}$ and $Q(\Acal)_{/X}$ to $Q^{inj}(\Acal)_{/X}$, so that it descends to $\Ncal$ by virtue of (i). Hence, $\Ncal^{inj}\to\Ncal$ is a homotopy equivalence. To conclude, remark that $\Ncal^{inj}$ has an initial object, given by the span $0\twoheadleftarrow 0\hookrightarrow X$.
\end{proof}
\begin{rmq} \label{DualResolution}
    The situation sometimes calls for a dual version of the above. Indeed, recall that $\Kth(\Ccal)\simeq\Kth(\catop{\Ccal})$ so that the Lemma also applies if $i$ is op-resolving, i.e. if the following two conditions are met:
    \begin{itemize}
        \item[(i')] For every exact sequence $\begin{tikzcd}[cramped]X\arrow[r, hook] & Y\arrow[r, two heads] & Z\end{tikzcd}$ with $Y\in\Acal$, we have $Z\in\Acal$
        \item[(ii)] For every $X\in\Ccal$, there is an exact sequence $\begin{tikzcd}[cramped]X\arrow[r, hook] & Y\arrow[r, two heads] & Z\end{tikzcd}$ with $Y\in\Acal$.
    \end{itemize}
\end{rmq}

The above result provides a counterpart to Theorem 3 of \cite{Quillen}. We do not prove, though it is true and the same proof suitably adapted to the realm of higher categories ought to apply, a version of Corollary 1. The reason for this is that in our setting, the situations we will consider will feature both resolving and op-resolving situation together. Hence, a corollary dedicated to composite of resolving functors makes little sense. 

However, some of Quillen's original arguments for Corollary 1 (the unnamed lemma that follows it) still manage to make their way into this paper, notably in the proof of Lemma \ref{QuillenTheoremThreeImpliesSameStab}. Hence, we are led to considering the aforementioned lemma and Theorem \ref{StableEnvelopeOfHeart} as the correct counterparts of this result for higher categories.

\section{Heart structures on stable categories}
\subsection{Generalities on heart structures}

\hspace{1.2em} Recall the definition of a weight structure on a stable category, owed to Bondarko for triangulated categories, our version for stable categories being inspired from the presentation in \cite{ElmantoSosnilo, HebestreitSteimleAppHarpaz}.

\begin{defi} \label{DefWeightStruct}
    Let $\Ccal$ be a stable category. A \textit{weight structure} on $\Ccal$ is the datum of a pair of full subcategories $(\Ccal_{\geq 0}, \Ccal_{\leq 0})$ subject to the following conditions:
    \begin{itemize}
        \item{(i)} $\Ccal_{\geq 0}, \Ccal_{\leq 0}$ are closed under retracts, $\Ccal_{\geq 0}$ under finite colimits and $\Ccal_{\leq 0}$ under finite limits.
        \item{(ii)} For any $X\in\Ccal_{\leq 0}$ and $Y\in\Ccal_{\geq 0}$, the mapping spectra $\map_\Ccal(X, Y)$ is connective
        \item{(iii)} For any $C\in\Ccal$, there is an exact sequence 
        $$
            \begin{tikzcd}[cramped]
                X\arrow[r] & C\arrow[r] & \Sigma Y
            \end{tikzcd}
        $$
        with $X\in\Ccal_{\leq 0}$ and $Y\in\Ccal_{\geq 0}$.
    \end{itemize}
    
    A \textit{weight-exact} functor between weight structures is an exact functor $f:\Ccal\to\Dcal$ on the underlying categories which preserves both non-negative and non-positively weighted objects. We denote $\CatExWT$ the categories of weight structures and weight-exact functors between them.
\end{defi}

If $(\Ccal_{\geq 0}, \Ccal_{\leq 0})$ is a weight structure on $\Ccal$ and $n,m\in\Z$, then, we let respectively $\Ccal_{\geq n}$, $\Ccal_{\leq m}$ and $\Ccal_{[m, n]}$ be the full subcategories of $\Ccal$ spanned by $\Sigma^n Y$ for $Y\in\Ccal_{\geq 0}$, by $\Omega^m X$ for $X\in\Ccal_{\leq 0}$ and their intersection. The category $\Ccal^{\heart}:=\Ccal_{[0,0]}$ is called the \textit{heart} of the weight structure. 

If $\Ccal$ is a weight structure, then, both $\Ccal_{\geq 0}$ and $\Ccal_{\leq 0}$ are closed under extensions as shown in Lemma 3.1.2 of \cite{HebestreitSteimleAppHarpaz}. This extends to the above defined categories.

\begin{defi} \label{BoundedDef}
    A weight structure on a stable category $\Ccal$ is said to be \textit{bounded} or \textit{exhaustive} if the map
    $$
        \bigcup_{n\in\Z} 
        \begin{tikzcd}
            \Ccal_{[-n,n]}\arrow[r, "\simeq"] & \Ccal
        \end{tikzcd}
    $$
    is an equivalence.
\end{defi}

We now introduce a generalization of Definition \ref{DefWeightStruct}: we remove the condition (ii) and replace the closure under retract by closure under extensions. The motivation for this definition is to encapsulate exactly what we will need for our proof of the theorem of the heart to work.

\begin{defi} \label{DefHeartStruct}
    Let $\Ccal$ be a stable category and $n\in\Z$. A \textit{heart structure} on $\Ccal$ is the datum of a pair of full subcategories $(\Ccal_{\geq 0}, \Ccal_{\leq 0})$ subject to the following conditions:
    \begin{itemize}
        \item{(i)} $\Ccal_{\geq 0}, \Ccal_{\leq 0}$ are closed under extensions, $\Ccal_{\geq 0}$ under finite colimits and $\Ccal_{\leq 0}$ under finite limits.
        \item{(iii)\footnote{The name of this condition is intentional, we simply consider condition (ii) empty.}} For any $C\in\Ccal$, there is an exact sequence 
        $$
            \begin{tikzcd}[cramped]
                X\arrow[r] & C\arrow[r] & \Sigma Y
            \end{tikzcd}
        $$
        with $X\in\Ccal_{\leq 0}$ and $Y\in\Ccal_{\geq 0}$.
    \end{itemize}
    
    A \textit{heart-exact} functor between heart structures is an exact functor $f:\Ccal\to\Dcal$ on the underlying categories which preserves the two full subcategories of the structure. We denote $\CatExHeart$ the category of heart structures and heart functors between them. 
\end{defi}

In particular, weight structures are heart structures and a functor between weight-exact categories is heart-exact if and only if it is weight-exact, so that $\CatExWT$ embeds as a full subcategory of $\CatExHeart$. Just as we say weighted category instead of a category with a weight-structure, we will say that $\Ccal$ is a heart category if it is a category endowed with a heart structure. \\

If $\Ccal$ is a heart category, then we say that $\Ccal$ is bounded if the same condition as Definition \ref{BoundedDef} holds for $\Ccal$. Moreover, we can also consider the intersection $\Ccal_{\geq 0}\cap\Ccal_{\leq 0}$, which we again denote $\Ccal^\heart$ and call the \textit{heart} of $\Ccal$.

\begin{rmq} \label{HeartInExact}
    If $F:\Ccal\to\Dcal$ is a heart-exact functor between heart categories, it maps $\Ccal^\heart$ to $\Dcal^\heart$ hence induces an additive functor $F^\heart:\Ccal^\heart\to\Dcal^\heart$. Since the exact structure on the heart is inherited from the stable structure, $F^\heart$ preserves the exact structures (i.e. is an exact functor of \textit{exact} categories) hence the heart is a well-defined functor landing in $\ExactCat$, the category of exact categories.

    We will see later on that a heart-exact $\Ccal\to\Dcal$ is in fact exactly the same datum as an exact functor $\Ccal^\heart\to\Dcal^\heart$, vindicating the name.
\end{rmq}

Let us record the following fact, which holds for weight structures by \cite[Lemma 3.1.5]{HebestreitSteimleAppHarpaz} and whose proof is exactly the same:

\begin{lmm} \label{BoundedSequencesInHeartStr}
    Let $X\in\Ccal_{[a;b]}$ and let $\begin{tikzcd}[cramped]Y\arrow[r] & X\arrow[r] & \Sigma Z\end{tikzcd}$ be an exact sequence with $Y\in\Ccal_{]-\infty; c]}$ and $Z\in\Ccal_{[c; +\infty[}$. 
    
    Then, we have $Y\in\Ccal_{[a; c]}$ and $Z\in\Ccal_{[c; b]}$.
\end{lmm}
\begin{proof}
    This is the same proof as \cite[Lemma 3.1.5]{HebestreitSteimleAppHarpaz} which shows the corresponding fact for weight structures: consider the rotated exact sequences:
    $$
    \begin{tikzcd}[cramped]
        Z\arrow[r] & Y\arrow[r] & X && X\arrow[r] & \Sigma Z\arrow[r] & \Sigma Y
    \end{tikzcd}
    $$
    The Lemma follows from closure under extensions of the relevant categories, which holds by hypothesis for us, applied to the above sequences.
\end{proof}

\subsection{On the category of bounded heart categories}

\hspace{1.2em} In \cite[Corollary 3.4]{SosniloHeart}, Sosnilo shows that taking the heart is a fully-faithful functor taking the category $\CatExWTBounded$ of bounded weighted categories to the full subcategory of $\CatAdd$ of weakly-idempotent additive categories. The goal of this section is to provide a similar result to Sosnilo's for heart categories. Instead of additive categories, our heart takes $\CatExHeartBounded$, the category of bounded heart categories, to in $\ExactCat$, the category of exact categories. We will show this functor is fully-faithful and although we do not have a complete description of the essential image, we will state some of its properties. \\

Recall that we have a fully-faithful functor $\CatEx\to\ExactCat$. If $\Ecal$ is an exact category, we can consider $\Stab(\Ecal^{split})$, the stable envelope of the underlying additive category of $\Ecal$ (see \cite[Construction 2.16]{BarwickGlasmanAkhilNikolaus}). There is a fully-faithful functor $\Ecal\to\Stab(\Ecal^{split})$ which is in general only additive and not exact. 

In \cite[Definition 3.1]{Klemenc}, Klemenc builds a stable category $\Stab(\Ecal)$ by localizing $\Stab(\Ecal^{split})$. He shows the composite $i:\Ecal\to\Stab(\Ecal)$ satisfies the following:

\begin{prop}[Klemenc] \label{StabLeftAdjoint}
    Let $\Ecal$ be an exact category. Then $i:\Ecal\to\Stab(\Ecal)$ is fully-faithful, closed under extensions, preserves and reflects exact sequences and induces an equivalence for every stable $\Dcal$:
    $$
        \begin{tikzcd}
            i^*:\FunEx(\Stab(\Ecal), \Dcal)\arrow[r, "\simeq"] & \FunEx(\Ecal, \Dcal)
        \end{tikzcd}
    $$
    In consequence, $\Stab:\ExactCat\to\CatEx$ is left adjoint to the inclusion.
\end{prop}
\begin{proof}
    The functor $i$ is fully-faithful by Proposition 3.17 of \cite{Klemenc}, exact by Corollary 3.21, has the universal property by 3.22, and reflects exact sequences and closed under extensions by Proposition 3.25.
\end{proof}

The fact that $\Ecal\to\Stab(\Ecal)$ is closed under extensions and reflects equivalences means in our lingo that the exact structure of $\Ecal$ is inherited from $\Stab(\Ecal)$. If $\Ecal$ is an exact 1-category, then $\Stab(\Ecal)$ is the bounded derived category of $\Ecal$ (see \cite[Corollary 3.29]{Klemenc}). \\

We now prove a more concrete criterion under which a functor between exact categories induces an equivalence of stable envelopes. 

\begin{lmm} \label{QuillenTheoremThreeImpliesSameStab}
    Let $i:\Acal\to\Ccal$ be a resolving functor. Then, it induces an equivalence 
    $$
        \begin{tikzcd}
            i^*:\FunEx(\Ccal, \Dcal)\arrow[r] & \FunEx(\Acal, \Dcal)
        \end{tikzcd}
    $$
    for every stable $\Dcal$. In particular, $i$ induces an equivalence $\Stab(\Acal)\xrightarrow{\simeq}\Stab(\Ccal)$.
\end{lmm}
\begin{proof}
    The crux of the proof lies in the fact that those conditions force the values on $\Ccal$ from the values on $\Acal$ for functors which preserve exact sequences. \\

    Let $X\in\Ccal$, then we define $\Ex(X)$ to be the category of exact sequences
    $$
        \begin{tikzcd}
            A\arrow[r, hook] & B\arrow[r, two heads] & X
        \end{tikzcd}
    $$
    where $A, B\in\Acal$. Note that $\Ex$ is not \textit{a priori} a functor in $X$. However, denote $\Null(X)$ the category of sequences $\begin{tikzcd}A\arrow[r] & B\arrow[r] & X\end{tikzcd}$ with a given null homotopy of their composite, where we no longer require the maps to be exact inclusions or exact projections. Remark that $\Null(X)$ is clearly a functor in $X$ as it is the fiber of the map
    $$
        \begin{tikzcd}
            \Fun(\{\bullet_2\longleftarrow\bullet_0\longrightarrow\bullet_1\}, \Ccal_{/X})\arrow[r] & \Fun(\{\bullet_2\}, \Ccal_{/X})\simeq\Ccal_{/X}
        \end{tikzcd}
    $$
    We claim the inclusion $\Ex(X)\to\Null(X)$ is cofinal. By \cite[Theorem 4.1.3.1]{HTT}, we are reduced to checking that for every $B_0\to X$ with $B_0\in\Acal$, the category of factorizations
    $$
        \begin{tikzcd}
            B_0\arrow[r] & B\arrow[r, two heads] & X
        \end{tikzcd}
    $$
    is weakly contractible. This category is nonempty: indeed, by (ii) there is a map $B\twoheadrightarrow X$ with $B\in\Acal$ and the induced $B_0\oplus B\twoheadrightarrow X$ is an exact projection by Lemma \cite[Lemma 4.7]{BarwickHeart}. Moreover, this category also admits products, given by the pullback $B\times_X B'\twoheadrightarrow X$ equipped with its canonical map from $B_0$. One checks that $B\times_X B'\in\Acal$ by closure under extension of $\Acal$ in $\Ccal$ and that the map is indeed a projections since they are stable under pullback and composition. Consequently, the comma category in question is indeed weakly contractible and the map cofinal. \\
    
    In consequence, if $F:\Acal\to\Ccal$ is any functor, the following formula defines a functor $R(F)$ with source $\Ccal$ and target $\Dcal$:
    $$
        X\longmapsto\colim_{\Ex(X)}\left(\cofib(F(A)\to F(B))\right)\simeq\colim_{\Null(X)} \left(\cofib(F(A)\to F(B))\right)
    $$
    Since colimits are functorial, $R$ upgrades to a functor $\Fun(\Acal, \Dcal)\to\Fun(\Ccal, \Dcal)$ taking $F$ to $R(F)$. The canonical $F(B)\to\cofib(F(A)\to F(B))$ induces a natural map in $X$ as follows:
    $$
        \colim_{\Null(X)}F(B)\longrightarrow\colim_{\Null(X)} \left(\cofib(F(A)\to F(B))\right)
    $$
    If $X\in\Acal$, then $\Null(X)$ has a terminal object given by the sequence $\begin{tikzcd}[cramped]0\arrow[r, hook] & X\arrow[r, two heads] & X\end{tikzcd}$ hence both colimits evaluate to $F(X)$, which provides a natural equivalence $F\to i^*R(F)$. This equivalence is again natural in $F$. Moreover, suppose $G:\Ccal\to\Dcal$ is a functor, there is a map $\cofib(G(A)\to G(B))\to G(X)$ for every object of $\Null(X)$, hence we have a natural transformation $R(i^*G)\to G$, which is itself natural in $G$. \\

    Suppose further that $F:\Acal\to\Dcal$ is exact, then we claim that $R(F)$ is also exact. We first show that $F$ sends exact sequences of $\Ex(X)$ to exact sequences, and then we deal with the more general case. Let 
    $$
    \begin{tikzcd}[cramped]
        A\arrow[r, hook] & B\arrow[r, two heads] & X && A'\arrow[r, hook] & B'\arrow[r, two heads] & X
    \end{tikzcd}
    $$ 
    be exact sequences in $\Ex(X)$. We have a diagram with exact rows and columns:
    $$
        \begin{tikzcd}
            0\arrow[r]\arrow[d] & A\arrow[r, equal]\arrow[d, hook] & A\arrow[d, hook] \\
            A'\arrow[r, hook]\arrow[d, equal] & B'\times_{X}B\arrow[r, two heads]\arrow[d, two heads] & B\arrow[d, two heads] \\
            A'\arrow[r, hook] & B'\arrow[r, two heads] & X
        \end{tikzcd}
    $$
    By closure under extension of $\Acal$, we have that $B'\times_{X}B\in\Acal$. Since $F$ is exact in $\Acal$, it sends every sequence save for the bottom horizontal and the right vertical ones to exact sequences in $\Dcal$. Since $\Dcal$ is stable, taking the pushout $F(B')\coprod_{F(B'\times_{X}B)}F(B)$ completes the diagram where we applied to $F$ to have exact rows and columns\footnote{Here it is critical that $\Dcal$ is stable, otherwise we would know nothing of the fiber of say $F(B)\to F(B'\times_{X}B)$.}. We deduce from this the following equivalences:
    $$
        \cofib(F(A)\to F(B))\simeq F(B')\coprod_{F(B'\times_{X}B)}F(B)\simeq\cofib(F(A')\to F(B'))
    $$
    This shows that the functor $\Ex(X)\to\Ccal$ sending $(\begin{tikzcd}[cramped]A\arrow[r, hook] & B\arrow[r, two heads] & X\end{tikzcd})$ to $\cofib(F(A)\to F(B))$ inverts every arrow in $\Ex(X)$.
    
    Since $\Ex(X)$ is cofinal in $\Null(X)$ which has an initial object, $\Ex(X)$ is contractible; it follows that the colimit over $\Ex(X)$ of $\cofib(F(A)\to F(B))$ is constant. In particular, for every $(\begin{tikzcd}[cramped]A\arrow[r, hook] & B\arrow[r, two heads] & X\end{tikzcd})$, the canonical map is an equivalence
    $$
        \cofib(F(A)\to F(B))\longrightarrow \colim_{\Ex(X)}\left(\cofib(F(A')\to F(B'))\right)
    $$
    Consequently, using the natural equivalence $F\to i^*R(F)$ we constructed, we see that $R(F)$ sends the objects of $\Ex(X)$ to exact sequences. \\

    We now deal with the general case and show $R(F)$ sends all the exact sequences in $\Ccal$ to equivalences. Let $\begin{tikzcd}[cramped]X\arrow[r, hook] & Y\arrow[r, two heads] & Z\end{tikzcd}$ be an exact sequence of $\Ccal$ and let $\begin{tikzcd}[cramped]A\arrow[r, hook] & B\arrow[r, two heads] & Y\end{tikzcd}$ be an exact sequence with $A, B\in\Acal$ as provided by hypothesis (i) and (ii) of Definition \ref{ResolvingDef}. Then, we have a diagram with exact rows and columns:
    $$
        \begin{tikzcd}
            A\arrow[r, equal]\arrow[d] & A\arrow[r]\arrow[d, hook] & 0\arrow[d] \\
            X\times_{Z} B\arrow[r, hook]\arrow[d, two heads] & B\arrow[r,two heads]\arrow[d, two heads] & Z\arrow[d, equal] \\
            X\arrow[r, hook] & Y\arrow[r, two heads] & Z
        \end{tikzcd}
    $$
    where $X\times_{Z} B\in\Acal$ by (i). Applying $R(F)$ to this diagram, every sequence save for the bottom horizontal one is sent to an exact one in $\Dcal$, hence this is also the case for the bottom horizontal one, which shows the wanted statement. Consequently, $R(F)$ is exact as wanted. \\

    We have a well-defined functor $R:\FunEx(\Acal, \Dcal)\to\FunEx(\Ccal, \Dcal)$ coming with a natural equivalence $\id\to i^*\circ R$. Moreover, if $G:\Ccal\to\Dcal$ is exact then the natural transformation $R(i^*G)\to G$ we have constructed earlier is an equivalence; indeed, we have shown previously that this is the case on $\Acal$ and for every $X\in\Ccal$, there is an exact sequence
    $$
        \begin{tikzcd}[cramped]
            A\arrow[r, hook] & B\arrow[r, two heads] & X
        \end{tikzcd}
    $$
    such that $A, B\in\Acal$; the result now follows from the exactness of both $G$ and $R(i^*G)$. Hence $R$ is an equivalence with inverse $i^*$, which proves the wanted claim.
\end{proof}

Of course, if $i:\Acal\to\Ccal$ is op-resolving, applying the above to $\catop{i}$ and $\catop{\Dcal}$ for a stable $\Dcal$ immediately implies the dual version:

\begin{cor} \label{QuillenOpResolvingImpliesSameStab}
    Let $i:\Acal\subset\Ccal$ be op-resolving. Then, the following functor
    $$
        \begin{tikzcd}
            i^*:\FunEx(\Ccal, \Dcal)\arrow[r] & \FunEx(\Acal, \Dcal)
        \end{tikzcd}
    $$
    is an equivalence for every stable $\Dcal$. In particular, $i$ induces an equivalence $\Stab(\Acal)\simeq\Stab(\Ccal)$.
\end{cor}

If $\Ccal$ is a heart category, then the inclusion $\Ccal^\heart\to\Ccal$ induces a functor $\Stab(\Ccal^\heart)\to\Ccal$. If $\Ccal$ is bounded, we now show that $\Ccal^\heart\to\Ccal$ factors as the composition of resolving or op-resolving functors. The above criterion then implies that they have the same stable envelope, but $\Ccal$ was already stable so $\Stab(\Ccal^\heart)\simeq\Ccal$.

\begin{thm} \label{StableEnvelopeOfHeart}
    Let $\Ccal$ be a bounded heart category. Then, the inclusion $\Ccal^\heart\to\Ccal$ factors as the (possibly infinite) composition of resolving or op-resolving functors. Consequently, we have an equivalence
    $$
        \begin{tikzcd}
            \Stab(\Ccal^\heart)\arrow[r, "\simeq"] & \Ccal
        \end{tikzcd}
    $$
\end{thm}
\begin{proof}
    The result follows from the following two facts:
    \begin{itemize}
        \item[(a)] Let $i\leq j$, then $\Ccal_{[i, j]}\to\Ccal_{[i, j+1]}$ is resolving
        \item[(b)] Let $i\leq j$, then $\Ccal_{[i, j]}\to\Ccal_{[i-1, j]}$ is op-resolving
    \end{itemize}
    Since the heart structure on $\Ccal$ is bounded, this implies that the inclusion $\Ccal^\heart\to\Ccal$ is a composition of resolving and op-resolving functors as wanted. The proof of (b) will be dual to the proof of (a), hence let us only do the latter. We remark already that $\Ccal_{[i, j]}$ is closed under extensions in $\Ccal$ hence in all of its exact subcategories. \\

    Let $\begin{tikzcd}[cramped]X\arrow[r, hook] & Y\arrow[r, two heads] & Z\end{tikzcd}$ be an exact sequence such that $X, Z\in\Ccal_{[i, j+1]}$ and $Y\in\Ccal_{[i,j]}$. Then, Lemma \ref{BoundedSequencesInHeartStr} ensures that $X\in\Ccal_{[i,j]}$ which gives (i). 

    If $Z\in\Ccal_{[i, j+1]}$, then applying condition (iii) of Definition \ref{DefHeartStruct} to $\Omega^j Z$ implies there is an exact sequence
    $$
        \begin{tikzcd}[cramped]X\arrow[r, hook] & Y\arrow[r, two heads] & Z\end{tikzcd}
    $$
    such that $X\in\Ccal_{\geq j}$ and $Y\in\Ccal_{\leq j}$. It follows from Lemma \ref{BoundedSequencesInHeartStr} that $Y\in\Ccal_{[i, j]}$, which gives (ii).
\end{proof}

\begin{cor} \label{FullyFaithfulHeart}
    The functor $(-)^\heart:\CatExHeartBounded\to\ExactCat$ taking a bounded heart category to its heart is fully-faithful.
\end{cor}
\begin{proof}
    Let $\Ccal, \Dcal$ be heart categories, we want to prove that the functor
    $$
        \begin{tikzcd}[cramped]
        \Phi:\Fun^{\heart-\Ex}(\Ccal, \Dcal)\arrow[r] & \FunEx(\Ccal^\heart, \Dcal^\heart)
        \end{tikzcd}
    $$
    is an equivalence, the left hand side denoting the full subcategory of $\FunEx(\Ccal, \Dcal)$ spanned by heart-exact functors. By Theorem \ref{StableEnvelopeOfHeart}, the inclusion $i:\Ccal^\heart\to\Ccal$ induces an equivalence:
    $$
        \begin{tikzcd}[cramped]
            i^*:\FunEx(\Ccal, \Dcal)\arrow[r, "\simeq"] & \FunEx(\Ccal^\heart, \Dcal)
        \end{tikzcd}
    $$
    Under this equivalence, the full subcategory of heart-exact functors is mapped to a full subcategory of $\FunEx(\Ccal^\heart, \Dcal^\heart)$. One checks that by the explicit formula for the inverse provided in Lemma \ref{QuillenTheoremThreeImpliesSameStab} that any functor $\Ccal^\heart\to\Dcal^\heart$ induces a heart-exact functor $\Ccal\to\Dcal$. This concludes.
\end{proof}

\begin{rmq}
    Differently stated, the above proposition says that an exact functor $\Ccal\to\Dcal$ between bounded heart categories is heart-exact if and only if it maps $\Ccal^\heart$ to $\Dcal^\heart$.
\end{rmq}

\subsection{Stably-exact categories}

\hspace{1.2em} We have seen in Theorem \ref{StableEnvelopeOfHeart} that the heart is fully-faithful. We are left with identifying the essential image.

\begin{defi}
    An exact category $\Ecal$ is said to be \textit{stably-exact} if it is the heart of a heart structure on a stable category.
\end{defi}

Recall that for any exact $\Ecal$, the exact structure on $\Ecal$ is inherited from $\Stab(\Ecal)$ so it makes sense to ask whether $\Ecal$ can be realized as an exact category via the heart of a heart category. Moreover, if $\Ecal$ is stably-exact, then it is the heart of a \textit{unique} heart structure on its stable envelope $\Stab(\Ecal)$ by Theorem \ref{StableEnvelopeOfHeart}.  \\

In the rest of this section, we provide examples and stability properties of the class of stably-exact categories.

\begin{prop}[Bondarko-Sosnilo] \label{BondarkoSosniloArgument}
    Additive categories which are closed under retracts in their stable envelopes, are stably-exact.
\end{prop} 
\begin{proof}
    This is essentially an argument that goes back to Bondarko. We spell it out here for two purposes: to translate it from the language of triangulated categories to stable categories, and to stress which part of the argument breaks when the exact structure is not split.
    
    Suppose $\Acal$ is a split-exact category. We write $\Ccal$ for $\Stab(\Acal)$ and we claim there is a heart structure on $\Ccal$ such that $\Ccal^\heart\simeq\Acal$. \\

    Denote $\Ccal_{\geq 0}$ the full subcategory of $\Ccal$ generated under finite colimits and extensions by objects $X\in\Acal$ (but not necessarily stable under limits), and dually for $\Ccal_{\leq 0}$ and limits. We check that this defines a bounded heart structure by proving the two points of Definition \ref{DefHeartStruct} are satisfied; remark that it is clear that such a structure is bounded. Point (i) is true by construction. We remark also that if $X, Y\in\Acal$ then $\map_\Ccal(X, Y)$ is connective; this clearly extends to $X\in\Ccal_{\leq 0}$ and $Y\in\Ccal_{\geq 0}$. \\

    \underline{(iii)}. We prove an apparently stronger result. Denote $\Wcal$ the full subcategory of $X\in\Ccal$ such that for every $n$, there exists an exact sequence
    $$
        \begin{tikzcd}
            P\arrow[r] & N\arrow[r] & X 
        \end{tikzcd}
    $$
    with $N\in\Ccal_{\leq n}$ and $P\in\Ccal_{\geq n}$; such a sequence will be called a $n$-shifted weight decomposition. We claim first that $\Acal\subset\Wcal$: indeed, for $n\geq 0$ we can take the trivial decomposition with $P=0$ and if $n\neq 0$, the following exact sequence is a $(-n)$-shifted weight decomposition
    $$
        \begin{tikzcd}
            \Omega^{n}X\oplus \Omega X\arrow[r] & \Omega^n X\arrow[r] & X
        \end{tikzcd}
    $$
    where the first map is the projection and the second map is the zero map. 

    To conclude, it remains to show that $\Wcal$ is closed under fibers and cofibers. The arguments are dual and so we consider the case of cofibers. Let $X\to Y$ be a map in $\Wcal$ and denote $Z$ its cofiber. Let $P\to N\to X$ be a $n$-shifted weight decomposition and $R\to Q\to Y$ a $(n+1)$-shifted weight decomposition. Remark that the composite $N\to Y$ factors uniquely through $Q$ since
    $$
        \begin{tikzcd}
            \map_\Ccal(N, Q)\arrow[r] & \map_\Ccal(N, Y)
        \end{tikzcd}
    $$
    has fiber $\map_\Ccal(N, R)$ where $N\in\Ccal_{\leq n}$ and $R\in\Ccal_{\geq n+1}$, and thus is an equivalence on $\pi_0$. Hence, by taking cofibers, we can form the diagram with exact rows and columns:
    $$
        \begin{tikzcd}
            P\arrow[r]\arrow[d] & R\arrow[r]\arrow[d] & T\arrow[d] \\
            N\arrow[r]\arrow[d] & Q\arrow[r]\arrow[d] & S\arrow[d] \\
            X\arrow[r] & Y\arrow[r] & Z
        \end{tikzcd}
    $$
    where $S\in\Ccal_{\leq n+2}$ and $T\in\Ccal_{\geq n}$. But remark that $T$ is an extension $R\to T\to\Sigma P$ of objects in $\Ccal_{\geq n+1}$, hence $T\in\Ccal_{\geq n+1}$ as well; similarly $S\in\Ccal_{\leq n+1}$, so that we have provided a $(n+1)$-shifted weight decomposition for $Z$. This shows that our structure satifies axiom (iii) of Definition \ref{DefHeartStruct}. \\

    Finally, we have to show that $\Acal$ is the intersection of $\Ccal_{\geq 0}$ and $\Ccal_{\leq 0}$. Recall that the Yoneda embedding $j:\Ccal\to\Ind\Ccal\simeq\FunEx(\catop{\Ccal}, \Sp)$ carries $\Acal$ to the full subcategory of compact objects in $\Pcal_\Sigma(\Acal)=\Fun^{\oplus}(\catop{\Acal}, \Sp_{\geq 0})$. It follows that $j$ also maps $\Ccal_{\geq 0}$ to compact objects of $\Pcal_\Sigma(\Acal)$, as connective spectra are closed under colimits in $\Sp$, and compact objects under finite colimits.
    
    Now, if $X\in\Ccal_{\geq 0}\cap\Ccal_{\leq 0}$, then $j(X)\in\Pcal_\Sigma(\Acal)$ and for every $Y\in\Acal$, we have $\map_{\Ind\Ccal}(j(X), j(Y))$ is connective, but the category $\Pcal_\Sigma(\Acal)$ is the non-abelian derived category of $\Acal$, in particular it is generated by $\Acal$ under sifted colimits \cite[Proposition 5.5.8.15]{HTT} thus the above connectivity extends to the whole of $\Pcal_\Sigma(\Acal)$. Hence, $\map_{\Pcal_\Sigma(\Acal)}(j(X), -)$ is connective and preserves sifted colimits; this extends to the mapping space because $\OmegaInfty:\Sp_{\geq 0}\to\Spaces$ preserve sifted colimits. In particular, since $\Pcal_\Sigma(\Acal)$ is generated by $\Acal$ under sifted colimits, this implies that $\id_{j(X)}$ factors through an object of $\Acal$. By hypothesis, $\Acal$ is closed under retracts in $\Ccal$, hence we get that $j(X)\in\Acal$ which concludes.  
\end{proof}

\begin{rmq}
    For any exact category, the category $\Wcal$ defined in the proof always contains $\Ecal$, and it is even closed under suspensions and loops. Hence, the reason the proof fails for general exact categories is the closure under extensions of $\Wcal$.
\end{rmq}

Let us now mention some stability properties for heart categories and stably-exact categories. If $\Ccal$ is stable, then it is stably-exact as an exact category, and combined with the above result, this means the two extremal type of exact categories are stably-exact. We now seek exhibit non-stable non-split-exact categories which are nonetheless stably-exact. First, let us mention the following fact:

\begin{lmm}
    Suppose $\Ecal$ is a stably-exact category, then $\catop{\Ecal}$ is stably-exact.
\end{lmm}
\begin{proof}
    We have an equivalence $\Stab(\catop{\Ecal})\simeq\catop{\Stab(\Ecal)}$. We claim that inverting the role of $\Stab(\Ecal)_{\geq 0}$ and $\Stab(\Ecal)_{\leq 0}$ gives a structure whose heart is $\catop{\Ecal}$. The weight decomposition of $X\in\Stab(\Ecal)$ yields a weight decomposition for $\Omega X\in\catop{\Stab(\Ecal)}$ (the loop being the loop of the opposite category), which concludes since $\Omega$ is an equivalence.
\end{proof}

\begin{lmm}
    Suppose $\Acal$ is a stably-exact category with stable envelope $\Ccal$, and denote $\Ccal_{[n, m]}:=\Ccal_{\geq n}\cap\Ccal_{\leq m}$ for integers $n\leq 0\leq m$. Then, $\Ccal_{[n, m]}$ is stably-exact. 
\end{lmm}
\begin{proof}
    We claim that the pair of categories $(\Ccal_{\geq n}, \Ccal_{\leq m})$ has the wanted properties of Definition \ref{DefHeartStruct}. Condition (i) is clear and condition (iii) is trivially realized by virtue of the inclusions $\Ccal_{\geq 0}\subset\Ccal_{\geq n}$ and $\Ccal_{\leq 0}\subset\Ccal_{\leq m}$, which imply we can just take the existing decompositions.
\end{proof}

Remark that even if $\Acal$ is additive, the $\Ccal_{[n, m]}$ need not be in general: this furnishes many examples of stably-exact categories which are neither additive nor stable.

\begin{lmm} \label{HeartStructureOnFunCat}
    Let $\Ccal$ be a bounded weight structure, then $\Fun(\Delta^n, \Ccal)$ has a bounded heart structure with heart $\Fun(\Delta^n, \Ccal^\heart)$.
\end{lmm}
\begin{proof}
    Suppose given a map $X\to Y$ in $\Ccal$ and exact sequences $P\to N\to X$ and $P'\to N'\to Y$ with $P, P'\in\Ccal_{\geq 0}$ and $N, N'\in\Ccal_{\leq 0}$. Then, there is a map $N\to N'$ making the square commute since
    $$
        \begin{tikzcd}
            \map_\Ccal(N, N')\arrow[r] & \map_\Ccal(N, Y)
        \end{tikzcd}
    $$
    is surjective on $\pi_0$ (its fiber is $\map_\Ccal(N, P')$ which is connective). In particular, this commutative square induces a map of the fibers $P\to P'$.
    
    Let $\Fun(\Delta^n, \Ccal)_{\geq 0}:=\Fun(\Delta^n, \Ccal_{\geq 0})$ and $\Fun(\Delta^n, \Ccal)_{\leq 0}=\Fun(\Delta^n, \Ccal_{\leq 0})$. Since limits and colimits are computed pointwise in functor categories, these are closed under extensions and finite limits or finite colimits in $\Fun(\Delta^n, \Ccal)$. Moreover, the above argument shows that the pointwise decompositions of (iii) upgrade to a decomposition in $\Fun(\Delta^n, \Ccal)$. Finally, it is clear that the heart is exactly $\Fun(\Delta^n, \Ccal^\heart)$.
\end{proof}

Remark that in the above, it is paramount that $\Ccal$ is a weight category, and not simply a heart category, and that the heart structure on $\Fun(\Delta^n, \Ccal)$ verifies a weakened version of axiom (ii) of weight categories, where mapping spectra are only supposed to be $(-n)$-connective. In particular, one has to consider $\Fun(\Delta^n, \Ccal^\heart)$ with its non-split exact structure for the results of the following sections to apply. \\

Finally, we conclude this section by a more algebro-geometric example of a heart structure. For a quasi-compact quasi-separated scheme $X$, we denote $\Perf(X)$ the category of compact $\Ocal_X$-modules.

We call a scheme $X$ \textit{divisorial} if it has an ample family of line bundles as defined in \cite[Definition 2.1]{ThomasonTrobaugh} (see also \cite[Exposé II, 2.2.3]{SGA6}). For instance, quasi-projective varieties over a field are divisorial by Example 2.1.2 of \cite{ThomasonTrobaugh}. Proposition 2.3.1(d) together with Theorem 2.4.3 of \textit{loc. cit.} show that perfect complexes over a divisorial $X$ are equivalent to bounded \textit{strict-complexes}, i.e. finite complexes of finite locally free $\Ocal_X$-modules. We will denote $\Vect(X)$ the subcategory\footnote{This is in fact a 1-category.} of such complexes concentrated in degree 0.

\begin{prop}
    Let $X$ be a divisorial scheme. Then, the category $\Perf(X)$ has a heart structure whose heart is $\Vect(X)$.
\end{prop}
\begin{proof}
    Define $\Perf(X)_{\geq 0}$ to be the full subcategory of $\Perf(X)$ of connective perfect $\Ocal_X$-modules and $\Perf_{Y}(X)_{\leq 0}$ those that are of negative tor-amplitude in $\Perf(X)$. 
    
    Since the two conditions are local, it follow from the affine case that both categories are closed under extensions, the former under finite colimits and the latter under finite limits. \\

    Let $K\in\Perf(X)$, by \cite{ThomasonTrobaugh} $K$ is a finite complex of finite locally free $\Ocal_X$-modules. Denote $P:=\tau_{\geq 1}K$ the positive truncation of $K$ and $N:=\tau_{\leq 0}K$ the non-positive one, so that the following sequence $N\to K\to P$ is an exact sequence in $\Perf(X)$. The freeness implies that $P$ is concentrated in strictly positive degrees, hence is the suspension of a connective complex. Moreover $N$ is necessarily of negative tor-amplitude in $\Perf(X)$ by \cite[\href{https://stacks.math.columbia.edu/tag/08CI}{Tag 08CI}]{StacksProject}, using that locally free modules are flat. Hence the previous exact sequence provides the decomposition of (iii) of Definition \ref{DefHeartStruct}. 
    
    Remark that a connective perfect complex of negative tor-amplitude is exactly concentrated in degree 0 by \textit{loc. cit.}, hence $\Vect(X)$ is indeed the heart of this heart structure. This concludes.
\end{proof}

\begin{rmq}
    In particular, the special case of \cite[Exercise 5.7]{ThomasonTrobaugh} where $i:Y\to X$ is the identity, shows that $\Kth(\Perf(X))$ agrees with the K-theory of $\Kth(\Vect(X))$, when viewing the latter as an exact 1-category and the former as a Waldhausen 1-category (see Lemma 3.5 of \textit{loc. cit.}). This is implied by the above proposition and Theorem \ref{HeartThm} below; in fact, the proof strategies are similar and Thomason-Trobaugh use the Gillet-Waldhausen theorem instead of our theorem of the heart. Note also that the main theorem of \cite{HiranouchiMochizuki} implies that $\Stab(\Vect(X))\simeq\Perf(X)$ with its exact structure, which we recover by the above and Theorem \ref{StableEnvelopeOfHeart}.
\end{rmq}

\section{Theorems of the heart}
\subsection{Resolution and the theorem of the heart for heart categories}

\hspace{1.2em} By following Quillen's original arguments, we have shown in the first section that the Q-construction and consequently K-theory, satisfied the resolution theorem, i.e. sends resolving (and op-resolving) functors to equivalence. Theorem \ref{StableEnvelopeOfHeart} shows if $\Ccal$ has a heart structure, then the map $\Ccal^\heart\to\Ccal$ factors as a composition of resolving and op-resolving functors. The theorem of the heart for heart structures immediately follows:

\begin{thm}[Theorem of the heart] \label{HeartThm}
    Let $\Ccal$ be a bounded heart structure and denote $\Ccal^\heart$ its heart. Then, the map
    $$
        \begin{tikzcd}
            \Kth(\Ccal^\heart)\arrow[r, "\simeq"] &\Kth(\Ccal)
        \end{tikzcd}
    $$
    is an equivalence, where on the left hand side, we have taken K-theory of the exact category $\Ccal^\heart$.
\end{thm}
\begin{proof}
    This follows from Theorem \ref{StableEnvelopeOfHeart} combined with Theorem \ref{QuillenTheoremThree} using that $\Kth$ preserves filtered colimits.
\end{proof}

By the fully-faithfulness of the heart functor, the above can be reformulated as follows:

\begin{cor}
    Let $\Ecal$ be a stably-exact category, then, the map
    $$
        \begin{tikzcd}
            \Kth(\Ecal)\arrow[r, "\simeq"] &\Kth(\Stab(\Ecal))
        \end{tikzcd}
    $$
    is an equivalence. 
\end{cor}

In fact, we claim a more general phenomenon is happening: for stably-exact categories, invariance under passage to the stable envelope is equivalent to inverting resolving and op-resolving functors.

\begin{thm} \label{HeartIsResolution}
    Let $F:\ExactCat\to\Ecal$ be a functor preserving filtered colimits. The following are equivalent:
    \begin{itemize}
        \item[(i)] For every stably-exact $\Ecal$, $F$ sends the map $\Ecal\to\Stab(\Ecal)$ to an equivalence.
        \item[(ii)] For every bounded heart category $\Ccal$, $F$ sends the map $\Ccal^\heart\to\Ccal$ to an equivalence.
        \item[(iii)] $F$ sends every resolving and every op-resolving functor $\Acal\to\Dcal$ with $\Acal, \Dcal$ stably-exact to an equivalence.
    \end{itemize}
\end{thm}
\begin{proof}
    Remark that (i) and (ii) are equivalent by Theorem \ref{StableEnvelopeOfHeart}, which also shows that (iii) implies (ii) since $F$ preserves filtered colimits. But Lemma \ref{QuillenTheoremThreeImpliesSameStab} shows (iii) implies (ii), which concludes.
\end{proof}

\subsection{The theorem of the heart for laced K-theory}

\hspace{1.2em} In this section, we draw some consequences of the previous generalization of the theorem of the heart for K-theory of endomorphisms, even allowing coefficients in a suitable bimodule. Let us first recall what they are. \\

Let $\Ccal$ be a stable category, a $\Ccal$-bimodule $M$ is an object of the category $\FunEx(\Ccal, \Ind\Ccal)$. This category is also equivalent to both $\EndL(\Ind\Ccal)$, the category of colimit-preserving endofunctors of $\Ind\Ccal$, and $\FunEx(\catop{\Ccal}\otimes\Ccal, \Sp)$, where  $\catop{\Ccal}\times\Ccal\to\catop{\Ccal}\otimes\Ccal$ is by definition the initial functor from $\catop{\Ccal}\times\Ccal$ which is exact in both variables. 

Every stable category admits a canonical nonzero bimodule, given by the Yoneda embedding (equivalently, $\id_{\Ind\Ccal}$ or $\map_\Ccal(-, -)$ in the two other pictures). 

\begin{defi}
    Let $M$ be a $\Ccal$-bimodule, we let $\Lace(\Ccal, M)$ denote the lax-equalizer of the Yoneda embedding $j$ to $M$, i.e. the following pullback of categories:
    $$
    \begin{tikzcd}
        \Lace(\Ccal, M)\arrow[r]\arrow[d] & \Ind(\Ccal)^{\Delta^1}\arrow[d] \\
        \Ccal\arrow[r, "{(j, M)}"] & \Ind(\Ccal)\times\Ind(\Ccal)
    \end{tikzcd}
    $$
    This is in particular a stable category by \cite[II.1.5]{NikolausScholze}.
\end{defi}
\begin{ex}
    Let $\Ccal:=\Perf(R)$ be the category of compact modules over a ring spectra. Then, a $\Perf(R)$-bimodule is fully-determined by its value on $R$, i.e. is of the form $M\otimes_R-$ where $M$ is a non-necessarily compact $R$-bimodule. \\

    Moreover, $\Lace(\Perf(R), M\otimes_R-)$ is a well-known object: it is exactly the category $\End(R, M)$ of $M$-parameterized endomorphisms whose objects are compact modules $N$ with a map $N\to M\otimes_RN$.
\end{ex}

The formalism of lax-equalizers \cite[Proposition II.1.5]{NikolausScholze} shows that the objects of $\Lace(\Ccal, M)$ are pairs $(X, f:X\to M(X))$ and that if $(X, f)$ and $(Y, g)$ are objects of $\Lace(\Ccal, M)$, the mapping spectra from $(X, f)$ to $(Y, g)$ in $\Lace(\Ccal, M)$ is given by the following equalizer:
$$
    \begin{tikzcd}[column sep=large]\map_\Ccal(X, Y)\arrow[r, shift left=1, "g_*"]\arrow[r, shift right=1, "f^*\circ M"'] & \map_{\Ind(\Ccal)}(X, M(Y))\end{tikzcd}
$$
Here we have suppressed every mention of the Yoneda embedding and view $\Ccal$ as a full subcategory of $\Ind\Ccal$. As a consequence of this formula, arrows in $\Lace(\Ccal, M)$ are the datum of an arrow in $\Ccal$ and a homotopy making the obvious square commute.

\begin{rmq}
    This lax-equalizer construction also features in \cite[Construction 3.2.5]{ElmantoSosnilo} for additive categories, though we warn the reader that in order to account for the shift in the equivalence $\Perf(R\oplus M)\simeq\Lace(\Perf(R), M\otimes_R-)$ for a connective ring spectra $R$ and a connective $R$-bimodule $M$ (see Example \ref{SqZeroVSLace} below), Elmanto and Sosnilo denote $\Ccal\oplus M$ the category $\Lace(\Ccal, \Sigma M)$.
\end{rmq}

We will show in an upcoming paper that $\Lace(\Ccal, M)$ deserves the name square-zero extension by identifying the tangent category of $\CatEx$ at some stable $\Ccal$ with the category of $\Ccal$-bimodules, and further showing that the canonical functor $\sqz:\TCatEx\to\CatEx$ is given by $\Lace$. \\

For the time being, we do not need this full structure and will be content by working over a fixed $\Ccal$. We will however need the following result:

\begin{lmm} \label{CoreLaceFibered}
    Let $M$ be a $\Ccal$-bimodule. The canonical functor $\core\Lace(\Ccal, M)\to\core\Ccal$ is the unstraightening of the functor $\core\Ccal\to\Spaces$ sending $X$ to $\Map(X, M(X))$.
\end{lmm}
\begin{proof}
    Consider $\Pcal$ the category given by the pullback
    $$
        \begin{tikzcd}[column sep=large]
            \Pcal\arrow[r]\arrow[d] & \TwAr(\Ind\Ccal)\arrow[d] \\
            \catop{\Ccal}\times\Ccal\arrow[r, "{(\catop{j}, M)}"] & \catop{\Ind(\Ccal)}\times\Ind\Ccal
        \end{tikzcd}
    $$
    For a general category $\Dcal$, $\TwAr(\Dcal)\to\catop{\Dcal}\times\Dcal$ classifies $\Map_{\Dcal}(-,-)$, hence it follows that $\Pcal\to\catop{\Ccal}\times\Ccal$ classifies the functor $(X, Y)\mapsto\Map(X, M(Y))$. Then, the lemma follows from the following pullback square:
    $$
        \begin{tikzcd}
            \core\Lace(\Ccal, M)\arrow[r]\arrow[d] & \Pcal\arrow[d] \\
            \core\Ccal\arrow[r, "\Delta"] & \catop{\Ccal}\times\Ccal
        \end{tikzcd}
    $$
    where $\Delta:\core\Ccal\to\catop{\Ccal}\times\Ccal$ is the diagonal $X\mapsto(X, X)$, using the canonical identification $\core\Ccal\simeq\core(\catop{\Ccal})$. By pasting, the above square is cartesian as soon as the following square also is:
    $$
        \begin{tikzcd}
            \core\Lace(\Ccal, M)\arrow[r]\arrow[d] & \TwAr(\Ind\Ccal)\arrow[d] \\
            \core\Ccal\arrow[r] & \catop{\Ind(\Ccal)}\times\Ind\Ccal
        \end{tikzcd}
    $$
    By the explicit description of pullbacks in $\CatEx$, we can replace the right vertical map by $\core\TwAr(\Ind\Ccal)\to\core\catop{\Ind(\Ccal)}\times\core\Ind\Ccal$. Now, the claim follows from the fact that $\core$ preserves pullbacks, that $\core\catop{\Ind(\Ccal)}\simeq \core\Ind(\Ccal)$ and $\core\TwAr(\Ind\Ccal)\simeq\core\Ar(\Ind\Ccal)$.
\end{proof}

If $\Ccal$ admits a heart structure, there is a subclass of bimodules $M$ which are suitably adapted to the weight structure. This is the following definition: 

\begin{defi}
    Let $\Ccal$ be a heart category. A $\Ccal$-bimodule $M$ is said to be a \textit{weighted bimodule} if for every $X\in\Ccal_{\leq 0}$ and $Y\in\Ccal_{\geq 0}$, the following mapping spectra
    $$
        \map_{\Ind\Ccal}(X, M(Y))
    $$
    is connective.
\end{defi} 
\begin{ex}
    Let $R$ be a connective ring spectra. Recall $\Perf(R)$ admits a weight structure, as it is the stable envelope of the additive $\Proj(R)$. A $R$-bimodule is weighted if and only for every connective $N$ and coconnective $P$, the spectrum
    $$
        \map_{\Mod_R}(P, M\otimes_R N)
    $$
    is connective. Applied to $N=R$, one checks that $M$ is necessarily connective, and this is also sufficient by the combination of the following two facts: whenever $P$ is coconnective and $Q$ is connective, $\map(P, Q)$ is itself connective and $M\otimes_R-$ preserves connectivity when $M$ is connective. Hence, weighted $R$-bimodules are exactly connective $R$-bimodules.
\end{ex}

In the situation where $\Ccal$ admits a weight structure and $M$ is a weighted bimodule, $\Lace(\Ccal, M)$ inherits a heart structure whose heart is $\Lace(\Ccal^\heart, M)$, the full subcategory of pairs $(X, f)$ with $X\in\Ccal^\heart$. This is not quite a weight structure, because the mapping spectra in $\Lace(\Ccal^\heart, M)$ need only be $(-1)$-connective; in fact, it is this fact that started the investigation which led to this article, and the introduction of heart structures.

\begin{lmm} \label{WeightStructureOnLace}
    Let $\Ccal$ be a heart category and $M$ a weighted bimodule. Then, $\Lace(\Ccal, M)$ admits a heart structure given by the full subcategories $\Lace(\Ccal_{\geq 0}, M)$ and $\Lace(\Ccal_{\leq 0}, M)$ fibered respectively over $\Ccal_{\geq 0}$ and $\Ccal_{\leq 0}$
\end{lmm}
\begin{proof}
    Clearly, $\Lace(\Ccal_{\geq 0}, M)$ and $\Lace(\Ccal_{\leq 0}, M)$ are closed under retracts and extensions in $\Lace(\Ccal, M)$, the former under pushouts and the latter under pullbacks.

    If $(Z, h:Z\to M(Z))$ is an object $\Lace(\Ccal, M)$, then there exists a weight decomposition $X\to Z\to \Sigma Y$, i.e. an exact sequence of $\Ccal$ with $X\in\Ccal_{\leq 0}$ and $Y\in\Ccal_{\geq 0}$. We have two diagrams
    $$
        \begin{tikzcd}
            Z\arrow[r, "h"]\arrow[d, "p"] & M(Z)\arrow[d, "M(p)"] & & X\arrow[r, dotted]\arrow[d, "i"] & M(X)\arrow[d, "M(i)"] \\
            \Sigma Y\arrow[r, dotted] & \Sigma M(Y) & & Z\arrow[r, "h"] & M(Z)
        \end{tikzcd}
    $$
    which it suffices to fill via the dotted arrows to get (iii) of Definition \ref{DefHeartStruct}. For the left hand side, we ought to show that the image of the map
    $$
        \map(\Sigma Y, \Sigma M(Y))\xrightarrow{p^*}\map(Z, \Sigma M(Y))
    $$
    contains $M(p)\circ h$. The cofiber of this map is exactly $\map(X, \Sigma M(Y))$ which is $1$-connective. Hence the above map is essentially surjective on $\pi_0$ providing the wanted dotted arrow. A dual argument deals with the other square. 
\end{proof}

\begin{rmq} \label{RemarkShiftedWeightOnLace}
    In the situation of the previous lemma, if the heart structure on $\Ccal$ is in fact a weight structure, $\Lace(\Ccal, M)$ verifies the following weaker verison of axiom (ii) of a weight structure:
    
    If $(X, f:X\to M(X))$ and $(Y, g:Y\to M(Y))$ with $X\in\Ccal_{\leq 0}$ and $Y\in\Ccal_{\geq 0}$, then
    $$
        \Map_{\Lace(\Ccal, M)}((X, f), (Y, g))
    $$ 
    is $(-1)$-connective. This follows from the aforementioned formula for this mapping spectra.
\end{rmq}

\begin{ex} \label{SqZeroVSLace}
    If $R$ is a connective ring spectrum and $M$ a connective bimodule, then $M\otimes_R-$ is weighted. Consequently, Remark \ref{RemarkShiftedWeightOnLace} shows that $\Lace(\Perf(R), \Sigma M\otimes_{R}-)$ has indeed a weight structure, whose heart is composed of objects of the form $(N, 0:N\to\Sigma M\otimes_{R}N)$ where $N\in\Proj(R)$. In particular, this heart is generated by one object as an additive category, namely $(R, 0:R\to \Sigma M)$. By Theorem \ref{StableEnvelopeOfHeart}, this implies that $\Lace(\Perf(R), \Sigma M\otimes_{R}-)$ is equivalent to $\Perf(S)$, where $S$ is the endomorphism ring spectrum of $(R, 0:R\to \Sigma M)$. 

    The underlying spectrum of $S$ is $R\oplus M$, as such an endomorphism is the datum of a $R$-linear map $R\to R$ and two nullhomotopies for $0:R\to\Sigma M$, i.e. a map $R\to\Omega\Sigma M=M$. It also follows that the induced ring structure is equivalent to the square-zero extension ring structure on $R\oplus M$. Hence, $\Lace(\Perf(R), \Sigma M\otimes_{R}-)\simeq\Perf(R\oplus M)$. Note that after applying K-theory, this equivalence was the starting point of Dundas-McCarthy's determination of stable K-theory in the celebrated \cite{DundasMcCarthy}.
\end{ex}

Applying Theorem \ref{HeartThm} to this situation, we have the following:

\begin{cor} \label{LacedHeart}
    Let $\Ccal$ be a bounded heart structure and $M$ a weighted bimodule. Denote $\Lace(\Ccal^\heart, M)$ the full subcategory of $\Lace(\Ccal, M)$ fibered over $\Ccal^\heart$. Then, there is an equivalence:
    $$
        \begin{tikzcd}
            \Kth(\Lace(\Ccal^\heart, M))\arrow[r, "\simeq"] & \Kth(\Lace(\Ccal, M))
        \end{tikzcd}
    $$
\end{cor}

In particular, when the heart structure on $\Ccal$ is a weight structure, the Yoneda embedding is one such weighted $\Ccal$-bimodule. Hence, the previous Corollary specifies to the following:

\begin{thm}[Theorem of the heart for $\Kth^{\End}$] \label{KendHeart}
    Let $\Ccal$ be a stable category equipped with a bounded weight structure. Then, there is an equivalence
    $$
        \begin{tikzcd}
            \Kth(\End(\Ccal^\heart))\arrow[r, "\simeq"] & \Kth(\End(\Ccal))
        \end{tikzcd}
    $$
\end{thm}

We finish this section by leveraging Corollary \ref{LacedHeart} to show that if $\Ccal$ has a bounded heart structure, then $\Kth(\Lace(\Ccal, M))$ commutes with sifted colimits taking values in weighted bimodules $M$. This kind of argument plays a key role in the proof of the Dundas-Goodwillie-McCarthy theorem (see \cite{DGMBook} for an account of this theorem in the setting of connective ring spectra), and will be used in a future article to provide a new proof of this theorem in the context of stable categories, without referring to ring spectra.

\begin{prop} \label{SiftedColimitMiracle}
    Let $\Ccal$ be a bounded heart category, then $\Kth(\Lace(\Ccal, M))$ preserves sifted colimits in $M$ which take values in weighted bimodules. 
\end{prop}
\begin{proof}
    Combining Lemma \ref{WeightStructureOnLace} and Theorem \ref{HeartThm}, we see that in the case where the sifted colimit has values in weighted bimodules, it suffices to prove that $\Kth(\Lace(\Ccal^\heart, M))$ commutes with sifted colimits. 
    
    Recall that $\LoopInfty:\Sp_{\geq 0}\to\Spaces$ preserves sifted colimits by \cite[1.4.3.9]{HA} hence it suffices to show the result for the space-valued delooped K-theory. Instead of the Q-construction, we will use the description of this space as the geometric realization of the core of the $S_\bullet$-construction, thanks to \cite[Theorem 3.10]{BarwickQConstr}. Since geometric realizations preserve colimits, it suffices to show that each $\core S_n\Lace(\Ccal^\heart, -)$ preserves sifted colimits of weighted bimodules. \\

    As a (necessary) warm-up, let us treat the first non-trivial case, $n=1$, where we ought to show that $\core\Lace(\Ccal^\heart, M)$ preserves sifted colimits in weighted $M$. We claim that for a fixed $X\in\Ccal^\heart$, the space $\Map_{\Ind\Ccal}(X, M(X))$ commutes with sifted colimits in weighted $M$. Indeed, this is clear at the level of the mapping spectra because $X$ is compact in $\Ind\Ccal$, and as we have already used $\LoopInfty:\Sp_{\geq 0}\to\Spaces$ preserves sifted colimits, so when $M$ is weighted, we have the first claim.
    
    But since $\core\Lace(\Ccal^\heart, M)$ is fibered over such spaces by \ref{CoreLaceFibered}, this extends to the whole space. Indeed, the unstraightening of a functor is its lax-colimit and lax-colimits commute with sifted colimits (combine the formula of \cite[Definition 2.9]{GepnerHausgengNikolaus} and the fact that sifted colimits of spaces commute to products). \\ 

    Recall that $S_n(\Lace(\Ccal^\heart, M))$ is the category of the following diagrams:
    $$
        \begin{tikzcd}
            0\arrow[r, hook] & (X_{1,1}, f_{1, 1})\arrow[r, hook]\arrow[d, two heads] & (X_{1,2}, f_{1,2})\arrow[r, hook]\arrow[d, two heads] & ...\arrow[r, hook]\arrow[d, two heads] & (X_{1, n}, f_{1,n})\arrow[d, two heads] \\
            & 0\arrow[r, hook] & (X_{2,2}, f_{2,2})\arrow[r, hook]\arrow[d, two heads] & ...\arrow[r, hook]\arrow[d, two heads] & (X_{2, n}, f_{2, n})\arrow[d, two heads] \\
            & & ...\arrow[r, hook] & ...\arrow[r, hook]\arrow[d, two heads] & ...\arrow[d, two heads] \\
            & & & 0\arrow[r, hook] & (X_{n, n}, f_{n,n})\arrow[d, two heads] \\
            & & & & 0
        \end{tikzcd}
    $$
    where $(X_{i, j}, f_{i,j})\in\Lace(\Ccal^\heart, M)$ and every square is exact. Of course, such a diagram is fully determined by its first row but remark that not all first rows need to induce a diagram where every object lies in $\Lace(\Ccal^\heart, M)$; in fact, this is already the case for $\Ccal^\heart$ in $\Ccal$.
    
    Denote $S_n(M)$ the induced bimodule on $S_n(\Ccal)$, obtained by applying $M$ to the diagrams of the above shape, then the above category is none other than $\Lace(S_n(\Ccal^\heart), S_n(M))$, viewed as a full subcategory of $\Lace(S_n(\Ccal), S_n(M))$. The arguments of the case $n=1$ will conclude, provided we can show that $S_n(M)$ is still a weighted bimodule, i.e. for every object $X:=(X_{i, j})\in\Lace(S_n(\Ccal^\heart), S_n(M))$ the mapping spectra $\map(X, S_n(M)(X))$ is connective. Such a mapping spectra is given by iterated pullbacks (because in $\Ccal$, the above diagrams are exactly equivalent to the datum of the first row), hence we have to show that
    $$
    \map(X_{1, 1}, M(X_{1, 1}))\times_{\map(X_{1, 1}, M(X_{1, 2}))}...\times_{\map(X_{1, n-1}, M(X_{1, n}))}\map(X_{1, n}, M(X_{1, n}))
    $$
    is connective, where all mapping spectra are taken in $\Ind\Ccal$. All of the terms appearing in the above are connective, so by induction, it suffices to show that every $\map(X_{1, i+1}, M(X_{1, i+1}))\to\map(X_{1, i}, M(X_{1, i+1}))$ is surjective on $\pi_0$. But $X_{1, i}\hookrightarrow X_{1, i+1}$ has its cofiber in $\Ccal^\heart$ and $M$ is exact and weighted, which concludes.
\end{proof}

\begin{rmq} \label{FailureHeartKEnd}
    Let us mention that we believe Corollary \ref{LacedHeart} (and thus subsequently Theorem \ref{KendHeart}) to be sharp insofar as if $\Ccal$ admits a bounded heart structure and $M$ is a bimodule which satisfies the following weaker version of being weighted:
    \begin{itemize}
        \item{(ii)} For any $X\in\Ccal_{\leq 0}$ and $Y\in\Ccal_{\geq 0}$, the mapping spectra $\map_{\Ind\Ccal}(X, M(Y))$ is $(-1)$-connective
    \end{itemize}
    then, we expect the statements to generally fail. Indeed, otherwise it would follow from considerations leading to the proof of the Dundas-Goodwillie-McCarthy theorem that cyclic K-theory, the fiber of $\Kth(\Lace(\Ccal; \Sigma M))\to\Kth(\Ccal)$, coincides with $\TR(\Ccal, \Sigma M)$, topological relative homology. In particular, it would follow that for $M$ a connective $R$-bimodule and $R$ a connective ring spectrum, $\Kcyc_0(R, M)\simeq\pi_0\TR(R, M)$; but as in explained in \cite{DottoKrauseNikolausPatchkoria}, $\pi_0\TR(R, M)$ identifies with the big Witt vectors whereas cyclic K-theory only identifies with \textit{rational} Witt vectors which in general only complete to the big Witt vectors for the $t$-adic topology, and otherwise do not coincide. We hope to prove the above considerations in an upcoming article.
\end{rmq}

\small
\bibliographystyle{alpha}
\bibliography{bibliographie}

\newcommand{\etalchar}[1]{$^{#1}$}
\begin{thebibliography}{CDH{\etalchar{+}}23d}

\bibitem[Bar13]{BarwickQConstr}
Clark Barwick.
\newblock On the {Q} construction for exact quasicategories.
\newblock {\em Available on the author's webpage.}, 2013.

\bibitem[Bar15]{BarwickHeart}
Clark Barwick.
\newblock On exact $\infty$-categories and the theorem of the heart.
\newblock {\em Compositio Mathematica}, 151(11):2160 --– 2186, 2015.

\bibitem[BGI71]{SGA6}
Pierre Berthelot, Alexander Grothendieck, and Luc Illusie.
\newblock {\em (SGA6) Théorie des {I}ntersections et {T}héorème de {R}iemann-{R}och}, volume 225 of {\em Lecture notes in mathematics}.
\newblock Springer-Verlag, 1971.

\bibitem[BGMN22]{BarwickGlasmanAkhilNikolaus}
Clark Barwick, Saul Glasman, Akhil Mathew, and Thomas Nikolaus.
\newblock {K}-theory and polynomial functors.
\newblock {\em Preprint}, 2022.

\bibitem[BGT13]{Blumberg}
Andrew~J. Blumberg, David Gepner, and Goncalo Tabuada.
\newblock A universal characterization of higher algebraic {K}-theory.
\newblock {\em Geometry \& Topology}, 17(2):733 -- 838, 2013.

\bibitem[BGT16]{BlumbergGepnerTabuadaKEnd}
Andrew~J. Blumberg, David Gepner, and Goncalo Tabuada.
\newblock {K}-theory of endomorphisms via noncommutative motives.
\newblock {\em Trans. Amer. Math. Soc.}, 368(2):1435–--1465, 2016.

\bibitem[BL23]{BurklundLevy}
Robert Burklund and Ishan Levy.
\newblock On the {K}-theory of regular coconnective rings.
\newblock {\em Selecta Mathematica}, 2023.

\bibitem[CDH{\etalchar{+}}23a]{HermKI}
Baptiste Calmès, Emanuele Dotto, Yonathan Harpaz, Fabian Hebestreit, Markus Land, Kristian Moi, Denis Nardin, Thomas Nikolaus, and Wolfgang Steimle.
\newblock Hermitian {K}-theory for stable $\infty$-categories {I}: {F}oundations.
\newblock {\em Selecta Mathematica}, 29(1):1 -- 269, 2023.

\bibitem[CDH{\etalchar{+}}23b]{HermKII}
Baptiste Calmès, Emanuele Dotto, Yonathan Harpaz, Fabian Hebestreit, Markus Land, Kristian Moi, Denis Nardin, Thomas Nikolaus, and Wolfgang Steimle.
\newblock Hermitian {K}-theory for stable $\infty$-categories {II}: {C}obordism categories and additivity.
\newblock {\em Preprint}, 2023.

\bibitem[CDH{\etalchar{+}}23c]{HermKIII}
Baptiste Calmès, Emanuele Dotto, Yonathan Harpaz, Fabian Hebestreit, Markus Land, Kristian Moi, Denis Nardin, Thomas Nikolaus, and Wolfgang Steimle.
\newblock Hermitian {K}-theory for stable $\infty$-categories {III}: {G}rothendieck-witt groups of rings.
\newblock {\em Preprint}, 2023.

\bibitem[CDH{\etalchar{+}}23d]{HermKIV}
Baptiste Calmès, Emanuele Dotto, Yonathan Harpaz, Fabian Hebestreit, Markus Land, Kristian Moi, Denis Nardin, Thomas Nikolaus, and Wolfgang Steimle.
\newblock Hermitian {K}-theory for stable $\infty$-categories {IV}: {P}oincaré motives and {K}aroubi-{G}rothendieck-{W}itt groups.
\newblock {\em Upcoming}, 2023.

\bibitem[DGM13]{DGMBook}
Bjorn Dundas, Tom Goodwillie, and Randy McCarthy.
\newblock {\em The local structure of algebraic K-theory}, volume~18.
\newblock Springer-Verlag London, 2013.
\newblock Algebra and Applications.

\bibitem[DKNP22]{DottoKrauseNikolausPatchkoria}
Emanuele Dotto, Achim Krause, Thomas Nikolaus, and Irakli Patchkoria.
\newblock Witt vectors with coefficients and characteristic polynomials over non-commutative rings.
\newblock {\em Compositio Mathematica}, 158(2):366--408, February 2022.

\bibitem[DM94]{DundasMcCarthy}
Bjørn~Ian Dundas and Randy McCarthy.
\newblock Stable {K}-theory and topological {H}ochschild homology.
\newblock {\em Annals of Mathematics}, 140(3):685--701, 1994.

\bibitem[Dun98]{DundasTheorems}
Bjørn~Ian Dundas.
\newblock K-theory theorems in topological cyclic homology.
\newblock {\em Journal of Pure and Applied Algebra}, 129(1):23--33, 1998.

\bibitem[ES21]{ElmantoSosnilo}
Elden Elmanto and Vladimir Sosnilo.
\newblock {On Nilpotent Extensions of $\infty$-Categories and the Cyclotomic Trace}.
\newblock {\em International Mathematics Research Notices}, 2022(21):16569--16633, 07 2021.

\bibitem[GHN17]{GepnerHausgengNikolaus}
David Gepner, Rune Hausgeng, and Thomas Nikolaus.
\newblock Lax colimits and free fibrations in $\infty$-categories.
\newblock {\em Documenta Mathematica}, 22:1225 -- 1266, 2017.

\bibitem[HHLN20]{HausgengHebestreitLinskensNuiten}
Rune Haugseng, Fabian Hebestreit, Sil Linskens, and Joost Nuiten.
\newblock Two-variable fibrations, factorisation systems and $\infty$-categories of spans.
\newblock {\em Preprint}, 2020.

\bibitem[HM10]{HiranouchiMochizuki}
Toshiro Hiranouchi and Satoshi Mochizuki.
\newblock {\em Pure weight perfect Modules on divisorial schemes}, pages 75--89.
\newblock Vieweg+Teubner, 2010.

\bibitem[HSH21]{HebestreitSteimleAppHarpaz}
Fabian Hebestreit, Wolfgang Steimle, and with an appendix by~Yonatan Harpaz.
\newblock Stable moduli space of hermitian forms.
\newblock {\em Preprint}, 2021.

\bibitem[Kle23]{Klemenc}
Jona Klemenc.
\newblock The stable hull of an exact $\infty$-category.
\newblock {\em Homology, Homotopy and Applications}, 24:195--220, 2023.

\bibitem[Lur08]{HTT}
Jacob Lurie.
\newblock {\em Higher Topos Theory}.
\newblock Princeton University Press, 2008.

\bibitem[Lur17]{HA}
Jacob Lurie.
\newblock {\em Higher Algebra}.
\newblock Princeton University Press, 2017.

\bibitem[Nee98]{NeemanHeartIIIA}
Anton Neeman.
\newblock {K}-theory for triangulated categories. {III(A). T}he theorem of the heart.
\newblock {\em Asian J. Math.2}, 3:495--–589, 1998.

\bibitem[Nee99]{NeemanHeartIIIB}
Anton Neeman.
\newblock {K}-theory for triangulated categories. {III(B). T}he theorem of the heart.
\newblock {\em Asian J. Math.2}, 3:557–--608, 1999.

\bibitem[NS17]{NikolausScholze}
Thomas Nikolaus and Peter Scholze.
\newblock On topological cyclic homology.
\newblock {\em Acta Mathematica}, 221, 2017.

\bibitem[Qui73]{Quillen}
Daniel Quillen.
\newblock Higher algebraic {K-Theory}: {I}.
\newblock {\em Algebraic {K}-theory (Proc. Conf., Northwestern Univ., Evanston, Ill., 1976)}, 1973.

\bibitem[Rap22]{Raptis}
George Raptis.
\newblock Dévissage for {W}aldhausen {K}-theory.
\newblock {\em Annals of K-theory}, 3:467--–506, 2022.

\bibitem[Sos19]{SosniloHeart}
Vladimir Sosnilo.
\newblock Theorem of the heart in negative k-theory for weight structures.
\newblock {\em Documenta Mathematica}, 24:2137 --– 2158, 2019.

\bibitem[{Sta}23]{StacksProject}
The {Stacks project authors}.
\newblock The stacks project.
\newblock \url{https://stacks.math.columbia.edu}, 2023.

\bibitem[TT90]{ThomasonTrobaugh}
Robert~W. Thomason and Thomas Trobaugh.
\newblock Higher algebraic {K-Theory} of schemes and of derived categories.
\newblock {\em The Grothendieck Festschrift, Vol. III}, 1990.

\bibitem[Wal85]{Waldhausen}
Friedhelm Waldhausen.
\newblock Algebraic {K-Theory} of spaces.
\newblock {\em Algebraic and Geometric Topology, Proceedings Rutgers 1983}, 1126:318 -- 419, 1985.

\end{thebibliography}

\end{document}